\documentclass[a4paper,12pt]{article}

\usepackage[pdftex]{graphicx}
\usepackage[section] {placeins}
\usepackage{enumerate}
\usepackage[normalem]{ulem}
\usepackage{amsmath}
\usepackage{amsthm}
\usepackage{latexsym}
\usepackage{amsfonts}
\usepackage{amssymb}

\usepackage{tikz-cd}
\usetikzlibrary{graphs}
\usetikzlibrary{arrows.meta}
\usetikzlibrary{shapes.geometric}

\usepackage{authblk}

\usepackage{latexsym,amssymb,amsfonts, amsthm, amsmath}
\usepackage[english]{babel}
\usepackage{bm}
\usepackage{mathrsfs}

\newtheorem{thm}{Theorem}[section]

\newtheorem{lem}[thm]{Lemma}
\newtheorem{cor}[thm]{Corollary}

\newtheorem{conj}[thm]{Conjecture}

\theoremstyle{definition} 
 
\newtheorem{exa}[thm]{Example}

\newcommand{\Z}{\mathbb{Z}}
\newcommand{\N}{\mathbb{N}}

\newcommand{\supp}{\mathop{\rm supp}\nolimits}

\newcommand*\boldell{\pmb{\ell}}

\newcommand{\apph}[3]{\mathrel{\operatorname*{{\uplus}^{#1}}_{#2}^{#3}}}

\usepackage{fullpage}

\begin{document}

\title{Construction Techniques for Linear Realizations of Multisets with Small Support }

\author[1]{Onur A\u{g}{\i}rseven \footnote{Corresponding author.  Email: \texttt{onura.marlboro.edu@gmail.com},  ORCID iD: 0009-0009-2056-1693}}
\author[2]{M.~A.~Ollis } 

\affil[1]{Unaffiliated, Boston, Massachusetts 02144, USA}
\affil[2]{Marlboro Institute for Liberal Arts and Interdisciplinary Studies, Emerson College, Boston, Massachusetts 02116, USA}

\maketitle

\begin{abstract}
A Hamiltonian path in the complete graph~$K_v$ whose vertices are labeled with the integers~$0,1,\ldots,v-1$ is a {\em linear realization} for the multiset $L$ of the linear edge-lengths (given by $|x-y|$ for the edge between vertices~$x$ and~$y$) of the edges in the path.  A linear realization is {\em standard} if an end-vertex is~0 and {\em perfect} if the end-vertices are~0 and~$v-1$.  

Linear realizations are useful in the study of the Buratti-Horak-Rosa (BHR) Conjecture on the existence of {\em cyclic realizations} (where cyclic edge-lengths are given by distance modulo $v$) for given multisets.  In this paper, we focus on multisets of the form~$\{1^a, (y-k)^b, y^c\}$.

Using core perfect linear realizations for supports of size 2 (which have the forms $\{x^{y-1},y^{x+1}\}$ whenever~$\gcd(x,y)=1$), we construct standard linear realizations  (with~$a=k-1$, $b=j(y-k)$, $c=jy$)  when~$k\mid y$ or $k \leq 4$.  When~$k=2$, these allow us to show that there is a linear realization whenever~$a \geq y$.  This is in line with the known results for the case of~$k=1$.  We also supplement these results for~$k=1$ by constructing linear realizations whenever~$b+c < y$ and $a \geq y - \min(b,c)$, from which the coprime version of the BHR Conjecture (requiring that~$v$ is coprime with each element of the multiset) follows for~$k=1$ when~$y \leq 16$.

Our methods show promise for constructing linear realizations for arbitrary $k$, in the direction of a resolution of the BHR Conjecture for supports of size 3.

\bigskip

\end{abstract}

\section{Introduction}\label{sec:intro}

Let~$K_v$ be the complete graph on~$v$ vertices and label the vertices~$\{0,1,\ldots v-1\} \subset \N$.   Define an induced edge-labelling, called the {\em linear length}, by labelling the edge between the vertices labelled~$x$ and~$y$ with~$\ell = |x-y| \in \N^{+}$.  Call an edge with label~$\ell$ an {\em $\ell$-edge}.

If~$L$ is the multiset of linear lengths in a subgraph~$\Gamma$ of~$K_v$, then~$\Gamma$ is a {\em linear realization} of~$L$.  A linear realization is {\em standard} if the first element is~0, while a standard linear realization is {\em perfect} if the final element has the largest label.  We are primarily interested in the case when~$\Gamma$ is a Hamiltonian path.

The main motivation for interest in linear realizations is their close relationship with cyclic realizations, which are obtained in the same way from a slightly different edge-labeling.  In~$K_v$, the {\em cyclic length} between the vertices~$x$ and~$y$ is given by distance modulo $v$; that is, $\min (|y-x| , v - |y-x|)$. This corresponds to viewing the labels as belonging to the additively written cyclic group~$\Z_v$. If the multiset of cyclic lengths in a subgraph~$\Gamma$ of~$K_v$ is~$L$ then~$\Gamma$ is a {\em cyclic realization} of~$L$.  We are again primarily interested in the case when~$\Gamma$ is a Hamiltonian path.

For~$1 \leq x < v$, let~$\widehat{x} = \min(x, v-x)$ be the {\em reduced form} of~$x$ with respect to~$v$.  A linear realization for the multiset~$\{ x_i : 1 \leq i \leq v-1 \}$ is a cyclic realization for the  multiset~$\{ \widehat{x_i} : 1 \leq i \leq v-1 \}$.  The two multisets are the same if and only if $x_i \leq  v/2$ for all~$i$.   Thus, although in the literature cyclic realizations are commonly called simply ``realizations", we will use the descriptors in this paper (except occasionally when we refer back to a specific realization or if the distinction is not necessary).

There are three main conjectures concerning cyclic realizations.  Note that the set $\supp(L) = \{ x : x \in L \}$ is the {\em support} of a multiset~$L$, and for each of the three conjectures, to say that it holds for a set~$U$ means that it holds for all multisets~$L$ with $\supp(L) = U$.

The first conjecture was proposed by Buratti in~2007.

\begin{conj}\label{conj:buratti}{\rm \cite{West07} (Buratti's Conjecture)}
Let~$v$ be prime.  If~$L$ is a multiset of size~$v-1$ with $\supp(L) \subseteq \{1, 2, \ldots,  \lfloor v/2 \rfloor\}$, then~$L$ is cyclically realizable.
\end{conj}

The second, proposed by Horak and Rosa and given this formulation by Pasotti and Pellegrini, extends Buratti's Conjecture to arbitrary values of~$v$.

\begin{conj}\label{conj:bhr}
{\rm \cite{HR09,PP14} (Buratti-Horak-Rosa (BHR) Conjecture)}
Let~$L$ be a multiset of size~$v-1$ with $\supp(L) \subseteq \{1, 2, \ldots,  \lfloor v/2 \rfloor\}$. Then~$L$ is cyclically realizable if and only if for any divisor~$d$ of~$v$ the number of multiples of~$d$ in~$L$ is at most~$v-d$.
\end{conj}

The third, proposed by the present authors, lies between the first two.  It considers those instances of the BHR Conjecture where the divisor condition is vacuous.

\begin{conj}\label{conj:cbhr}{\rm \cite{AO2} (Coprime BHR Conjecture)}
If~$L$ is a multiset of size~$v-1$ such that $\supp(L) \subseteq \{1, 2, \ldots,  \lfloor v/2 \rfloor\}$ and $\gcd(v,x) = 1$ for all~$x \in L$, then~$L$ is cyclically realizable.
\end{conj}

The BHR Conjecture, and hence the other two conjectures, have implications for graph decomposition problems, particularly in relation to the cyclic group~$\Z_v$, see~\cite{CD10,OPPS} for exploration of these implications.  The conjectures are known to hold for supports of size at most~2~\cite{DJ09,HR09}.  In this paper, we focus on multisets with supports of size~3. We prove new results for a variety of those that are, in the sense we describe below, among the furthest from being known to satisfy the Coprime BHR Conjecture.

A tool that is especially powerful in the context of the Coprime BHR Conjecture is that of equivalence.  Since each  automorphism of~$\Z_v$ is multiplication by some element~$s \in \Z_v$ with~$\gcd(s,v) = 1$, applying an automorphism sends a cyclic edge-label~$x$ to the cyclic edge-label~$\widehat{sx}$.  If we  take~$s = x^{-1}$ for some~$x$ in the multiset then we have~1 as an element of the resulting equivalent multiset.

Therefore, when considering the Coprime BHR Conjecture for supports of size~3, we may assume, without loss of generality, that our multiset is of the form $L = \{1^a, x^b, y^c\}$ with $x, y \in \N^{+}$ and $a, b, c \in \N$.  (We use exponents in multisets to indicate multiplicity.)  For this type of multiset, Theorem~\ref{th:1xy_ao2} provides a strong general result that depends only on~$a$ and the parities of~$x$ and~$y$. See~\cite{AO2} for a full description including further results that consider various special cases.

\begin{thm}\label{th:1xy_ao2}{\rm \cite{AO2,OPPS} (Select known results)}
Let $L = \{1^a, x^b, y^c\}$. Define a function
$$f(x,y) = \begin{cases}
  y-1 & \text{if $x$ and $y$ are even} \\
  x+y-2 & \text{if $x=3$ or if $x$ is even and $y$ is odd,} \\
  x+y-1 & \text{otherwise.}
\end{cases}$$
If~$a\geq f(x,y)$ then~$L$ is linearly realizable.  
\end{thm}

A distinctive, and useful, quality of Theorem~\ref{th:1xy_ao2} is that it is effective when~$a$ is sufficiently large, without further caveats.  Using the idea of equivalence, when considering the Coprime BHR Conjecture, we have three equivalent multisets that contain~1.  If the number of 1s in any one of these multisets is sufficiently large, in the sense given by Theorem~\ref{th:1xy_ao2}, then each is realizable.
Using this approach, we arrive at Theorem~\ref{th:1xy_new}.

\begin{thm}\label{th:1xy_new}{\rm \cite{AO2} (Counterexample characterization)}
Let $L = \{1^a, x^b, y^c\}$ be a multiset of size~$v-1$ with equivalent multisets $L' = \{1^b, \dot{x}^c, \dot{y}^a  \}$ and $L'' = \{1^c, \ddot{x}^a, \ddot{y}^b\}$.   Let~$f$ be defined as in Theorem~\ref{th:1xy_ao2}.

If~$a\geq f(x,y)$, $b \geq f(\dot{x}, \dot{y})$ or $c \geq f( \ddot{x} , \ddot{y})$ then~$L$ is realizable.  In particular, any counterexample to the Coprime BHR Conjecture has
$$ f(x , y) + f(\dot{x} ,\dot{y}) + f(\ddot{x} , \ddot{y}) \geq v+2. $$
\end{thm}

Sometimes, Theorems~\ref{th:1xy_ao2} and~\ref{th:1xy_new} are sufficient to prove the Coprime BHR Conjecture for a support of size~3.  More commonly, they restrict the ranges of outstanding cases.

\begin{exa}\label{ex:mot}{\rm (Uses and limitations of the characterization)}
Let us consider the support~$\{1,17,19\}$ with respect to two values of~$v$:~103 and~105.   

When~$v = 103$, the multiset $L = \{ 1^a, 17^b, 19^c \}$ is equivalent to $L' = \{1^b, 6^a, 11^c\}$, and $L'' = \{ 1^c , 28^a , 38^b \}$.  In the notation of Theorem~\ref{th:1xy_new} we find
$$ f(x , y) + f(\dot{x} ,\dot{y}) + f(\ddot{x} , \ddot{y}) = 35 + 15  + 37 = 87 < 105 = v+2. $$
Hence, for~$v = 103$, the BHR Conjecture holds for the support~$\{1,17,19\}$.

When~$v = 105$, the  multiset $L = \{ 1^a, 17^b, 19^c \}$ is equivalent to $L' = \{1^c, 11^a, 23^b\}$, and $L'' = \{ 1^b , 32^c , 37^a \}$.  In the notation of Theorem~\ref{th:1xy_new} we find
$$ f(x , y) + f(\dot{x} ,\dot{y}) + f(\ddot{x} , \ddot{y}) = 35 + 33  + 67 = 135 > 107 = v+2. $$
Hence, for~$v = 105$, Theorem~\ref{th:1xy_new} is not strong enough to allow us to conclude that the BHR Conjecture holds for the support~$\{1,17,19\}$. However, using Theorem~\ref{th:1xy_ao2}, we can deduce that any counterexample to the conjecture must have~$a \leq 34$, $b \leq 66$ and~$c \leq 32$.  Further, as $a+b+c = v-1$, we also know that a counterexample must have~$a = v-1-b-c \geq 105 - 1 - 66 - 32 = 6$.
\end{exa}

The second situation in Example~\ref{ex:mot} illustrates how we informally think about which supports of size~3 are furthest from being known to satisfy the conjecture.  Take~$L$ and $f$ as in Theorem~\ref{th:1xy_ao2}.  If~$L$ has $a < f(x,y)$ then we cannot use Theorem~\ref{th:1xy_ao2} directly.  However, the smaller that~$a$ is, the larger that~$b$ or~$c$ must be, and  the likelihood that we can apply Theorem~\ref{th:1xy_new} increases as~$b$ and~$c$ increase.  We therefore think of instances where~$a < f(x,y)$ with larger~$a$ as the most troublesome. 

The bound given by~$f$ is largest when at least one of~$x$ and~$y$ is odd and~$x$ is close to~$y$ in size, hence our focus on multisets with support~$\{1, y-k, y\}$ for small~$k$.  We take $0 \leq k < y$ throughout the paper.  Our methods apply also when~$x$ and~$y$ are both even, thus give an alternative approach to~\cite{OPPS}.

The case of $k=1$ was studied in~\cite{AO1}, where it was shown that for support~$\{1,y-1,y\}$ the Coprime BHR Conjecture holds (i) in all cases when~$v \geq 2y^2 + 9y$, (ii) in almost all cases when~$y \leq 16$.  In this paper, we complete the work on the cases when~$y \leq 16$, and introduce new techniques to extend our methods to the~$k=2$ case, paving the way for future investigations for general $k$.

Section~\ref{sec:background} reviews and expands the existing theory and constructions to be used.  Sections~\ref{sec:para} and~\ref{sec:1y-1y} together provide constructions for multiple new types of standard and perfect linear realizations.  Section~\ref{sec:1y-1y} revisits the case of $k=1$, completing the work on cases with small~$y$.  Section~\ref{sec:y-2} produces our main result for $k=2$. Section~\ref{sec:future} identifies the next steps needed for more general results for $k\geq 2$, akin to those from~\cite{AO1} for $k=1$.

Theorem~\ref{th:constr} below lists the essential new constructions (from Theorems~\ref{th:2perf} and \ref{th:tracks23}, and Lemma~\ref{lem:sawtooth}, respectively) leading to our main results. Note that the perfect linear realizations constructed for supports of size 2 are a considerable generalization of those currently known, and applicable much more widely than their uses in this paper.

\begin{thm}\label{th:constr}{\rm (Select new constructions)} Let $x,y>1$, $k \leq 4$, and $b,c<y-1$.

\begin{itemize}
\item $L = \{ x^{y-1}, y^{x+1} \}$ has a perfect linear realization if $x$ and~$y$ are coprime.  

\item $L = \{ 1^{k-1}, (y-k)^{j(y-k)} , y^{j y} \}$ has a standard linear realization for all~$j \geq 1$.

\item $L = \{1^a, (y-1)^b, y^c\}$ has a standard linear realization when~$a \geq y - \min(b,c)$.
\end{itemize}

\end{thm}

Theorem~\ref{th:small_y-1} is the completion of the earlier work on the $k=1$ case for~$y \leq 16$, while Theorem~\ref{th:biggie} provides the current form of our main result by combining Lemma~\ref{lem:omega2} ($k=0$), Theorem~\ref{th:y-1} ($k=1$), and Theorem~\ref{th:y-2} ($k=2$).

\begin{thm}\label{th:small_y-1}{\rm ($k=1$; small $y$)}
The Coprime BHR Conjecture holds for supports of the form~$\{1,y-1,y\}$ for all~$y \leq 16$.
\end{thm}

\begin{thm}\label{th:biggie}{\rm ($k \leq 2$; $a \geq y$)}
~$L = \{ 1^a , (y-k)^b, y^c \}$ has a standard linear realization for $k \leq 2$ whenever~$a \geq y$.
\end{thm}

The next example demonstrates how Theorem~\ref{th:biggie} helps prove the Coprime BHR Conjecture for more multisets of support of size~3, sometimes completing the proof for a support at a given $v$.

\begin{exa}\label{ex:mot2}{\rm (Some use cases of the main result)} Consider again the multiset $L = \{1^a, 17^b, 19^c\}$ from Example~\ref{ex:mot}.

When~$v = 105$, Theorems~\ref{th:1xy_ao2} and~\ref{th:1xy_new} had shown that any potential counterexample to the Coprime BHR Conjecture must have $6 \leq a \leq 34$.  Theorem~\ref{th:biggie} reduces this range to $6 \leq a \leq 18$.

When~$v = 127$, Theorems~\ref{th:1xy_ao2} and~\ref{th:1xy_new} can be shown to imply that a counterexample must have~$24 \leq a \leq 34$.  Theorem~\ref{th:biggie} covers this range.

\end{exa}

We believe the next logical extension of our methods could provide a resolution to the following conjecture, which would be a significant step towards a proof of the Coprime BHR Conjecture for supports of size~3.

\begin{conj}\label{conj:y-k}{\rm (Linear realizations for small $k$)}
When~$k < y/2$, multisets of the form~$\{1^a, (y-k)^b, y^c \}$ are linearly realizable for all~$a \geq y$.  
\end{conj}

\section{Notation, tools, and previous work}\label{sec:background}

In this section, we review and expand the existing theory and constructions that we require. We start with two basic manipulations of linear realizations that allow us to connect them together in two specific ways. These connecting methods, in return, lead to convenient constructions for supports of size 2, with useful properties of extendibility, as discussed toward the end of the section. We also provide visualization methods for these constructions.

For each linear realization, we have two particularly useful manipulations: Its translation and its complement.  The {\em (embedded) translation} by~$t$ is formed by adding~$t$ to each vertex.  It realizes the same multiset, as a non-Hamiltonian path in~$K_{v+t}$.  The {\em complement} of a linear realization is obtained by replacing each vertex~$x$ with~$v-1-x$.  It realizes the same multiset, as a Hamiltonian path in~$K_{v}$.  

An important advantage of working with linear realizations over cyclic ones is the possibility to connect them together without disturbing the internal lengths.  Two standard linear realizations can be connected together via either their concatenation or a bridge concatenation.

Concatenation is done by identifying a pair of end-vertices as follows:  Suppose~$L$ is a multiset of size~$v-1$ with standard linear realization~$\mathbf{g}$, and that~$M$ is a multiset of size~$w-1$ with standard linear realization~$\mathbf{h}$. 
Then the {\em concatenation} of~$\mathbf{g}$ and~$\mathbf{h}$, denoted~$\mathbf{g} \oplus \mathbf{h}$,  is the Hamiltonian path in~$K_{v+w-1}$ obtained by taking complement of~$\mathbf{g}$ and the translation of~$\mathbf{h}$ by~$v-1$ and identifying the vertices labeled~$v-1$. 
A concatenation has the properties described in Lemma~\ref{lem:concat}.

\begin{lem}\label{lem:concat}{\rm \cite{HR09,OPPS} (Concatenation)}
Suppose that~$L$ and~$M$ are multisets with standard linear realizations~$\mathbf{g}$ and~$\mathbf{h}$ respectively. Then~$\mathbf{g} \oplus \mathbf{h}$ is a linear realization of~$L \cup M$.  If~$\mathbf{h}$ is perfect, then~$\mathbf{g} \oplus \mathbf{h}$ is standard;  if both~$\mathbf{g}$ and~$\mathbf{h}$ are perfect, then~$\mathbf{g} \oplus \mathbf{h}$ is perfect.
\end{lem}

As a special case of concatenation, we have the prepending of 1-edges, which, in many situations, lets us move from a construction of a linear realization for a particular multiset $\{1^{a'}, x^b, y^c\}$ to a construction of a linear realization for  $\{1^{a}, x^b, y^c\}$ for each~$a \geq a'$.

\begin{lem}\label{lem:concat1}{\rm (Prepending 1-edges)}
If~$L$ has a standard  linear realization then for all~$s \geq 0$ the multiset~$L \cup \{1^s\}$ has a standard  linear realization.  If the linear realization for~$L$ is perfect, then so is the linear realization for~$L \cup \{1^s\}$. 
\end{lem}

\begin{proof}[Proof Construction.]
Let~$\mathbf{h}$ be the standard (or perfect) linear realization of~$L$ and concatenate it with~$\mathbf{g} = [ 0, 1, \ldots, s ] $, which is a  perfect linear realization for~$\{1^s\}$. 
\end{proof}

Bridge concatenation, on the other hand, is done via the addition of a new edge between two end-vertices as follows:  In~$K_v$, take two paths~$\mathbf{p_1}$ and~$\mathbf{p_2}$ with end-points~$g$ and~$h$ respectively and denote the linear length  of the edge between~$g$ and~$h$ by~$\ell$.  The path obtained by adding the $\ell$-edge between~$g$ and~$h$ is a {\em bridge concatenation}, denoted~$\mathbf{p_1} \uplus^\ell \mathbf{p_2}$.  This construction is not necessarily unique, but in all of the use cases we encounter in this paper, the bridging choice is clear from context.

More generally, given paths~$\mathbf{p_1}, \mathbf{p_2}, \ldots, \mathbf{p_t}$ and a vector of lengths~$\boldell  = (\ell_1, \ell_2, \ldots, \ell_{t-1})$,  let
$$ \apph{\boldell}{i=1}{t}  \mathbf{p_i}  = 
\mathbf{p_1} \uplus^{\ell_1} \mathbf{p_2} \uplus^{\ell_2} \cdots \uplus^{\ell_{t-1}} \mathbf{p_t}
{\rm \ \ \ \  and \ \ \ \  } 
\apph{\boldell}{i=t}{1}  \mathbf{p_i} = 
\mathbf{p_t} \uplus^{\ell_1} \mathbf{p_{t-1}} \uplus^{\ell_2} \cdots \uplus^{\ell_{t-1}} \mathbf{p_1}.$$

A {\em spanning linear forest of~$K_v$} is a collection of disjoint paths whose vertices partition the vertex set.
A typical use of the bridge concatenation construction is to take the set of paths $\mathbf{p_1}, \mathbf{p_2}, \ldots, \mathbf{p_t}$ to be a spanning linear forest of~$K_v$ chosen in such a way that the ends of the path may be connected by bridge concatenations using 1-edges.  

This approach was used in~\cite{AO1}, where the spanning linear forest consists of ``fauxsets".  Let~$v= qx+r$ with~$0 \leq r < q$.  Define the {\em fauxset} $\phi_i$ with respect to~$x$ to be the subset of vertices in~$K_v$ that are congruent to~$i \pmod{x}$.  Let the {\em fauxset traversal} $\Psi_i$ be the path that traverses the fauxset in sequence, realizing only $x$-edges.

It is sometimes useful to have concise notation working with paths in fauxsets.  Following~\cite{AO2}, let
$$q^* = q^*(v,r,i) =  
\begin{cases} 
q & \mathrm{if \  } i<r, \\
 q-1 & \mathrm{otherwise.}
\end{cases}
$$
That is, $q^*$ is the value that makes $q^*x + i$ the largest element of~$\varphi_i$.  Let the {\em fauxset segment} $\Psi_i^{(a,b)}$ be the path $[ ax + i, (a+1)x+i, \ldots, bx+i ]$ of elements in~$\varphi_i$ that realizes~$\{x^{b-a} \}$.  If the fauxset segment has the single vertex~$ax+i \in \phi_i$,  we denote it~$\Psi_i^{(a)}$.

We now have the tools and notation to introduce the $\omega$-constructions, a collection of linear realizations with support~$\{1,x\}$ constructed by using 1-edges to connect the spanning linear forest given by fauxsets with respect to~$x$.  Most $\omega$-constructions can be built using only the full fauxset traversals. For one specific type, it's necessary to traverse one fauxset via segments near one end of the Hamiltonian path.  These two methods are described in Lemmas~\ref{lem:omega1} and~\ref{lem:omega2}, respectively.

\begin{lem}\label{lem:omega1}{\rm \cite{AO1} ($\omega$-constructions)}
There is a standard linear realization for~$\{1^{x-1}, x^b\}$, except when, in the Euclidean division~$b = q'x + r'$ for~$b \geq 1$ with~$0 \leq r' < x$, both~$x$ and~$r'$ are odd with~$r'>1$.   
\end{lem}

\begin{figure}
\caption{Standard linear realizations~$\mathbf{h_1}$ and~$\mathbf{h_2}$ for~$\{1^{x-1}, x^b\}$ for $(v,x,b) = (25, 7,18)$, $(29, 8,21)$, $(12, 7, 5)$ and $(22, 7,15)$  {\rm \cite{AO1}}.}\label{fig:omega1}
\begin{center}
\begin{tikzpicture}[scale=0.9, every node/.style={transform shape}]
\fill (1,1) circle (2pt) ; \fill (1,2) circle (2pt) ; \fill (1,5) circle (2pt) ; \fill (1,6) circle (2pt) ; \fill (1,7) circle (2pt) ; \fill (1,8) circle (2pt) ; \fill (2,1) circle (2pt) ; \fill (2,2) circle (2pt) ; \fill (2,5) circle (2pt) ; \fill (2,6) circle (2pt) ; \fill (2,7) circle (2pt) ; \fill (2,8) circle (2pt) ; \fill (3,1) circle (2pt) ; \fill (3,2) circle (2pt) ; \fill (3,5) circle (2pt) ; \fill (3,6) circle (2pt) ; \fill (3,7) circle (2pt) ; \fill (3,8) circle (2pt) ; \fill (4,1) circle (2pt) ; \fill (4,2) circle (2pt) ; \fill (4,5) circle (2pt) ; \fill (4,6) circle (2pt) ; \fill (4,7) circle (2pt) ; \fill (4,8) circle (2pt) ; \fill (5,1) circle (2pt) ; \fill (5,5) circle (2pt) ; \fill (5,6) circle (2pt) ; \fill (5,7) circle (2pt) ; \fill (6,1) circle (2pt) ; \fill (6,5) circle (2pt) ; \fill (6,6) circle (2pt) ; \fill (6,7) circle (2pt) ; \fill (7,0) circle (2pt) ; \fill (7,1) circle (2pt) ; \fill (7,5) circle (2pt) ; \fill (7,6) circle (2pt) ; \fill (7,7) circle (2pt) ; \fill (9,1) circle (2pt) ; \fill (9,2) circle (2pt) ; \fill (9,3) circle (2pt) ; \fill (9,5) circle (2pt) ; \fill (9,6) circle (2pt) ; \fill (9,7) circle (2pt) ; \fill (9,8) circle (2pt) ; \fill (10,1) circle (2pt) ; \fill (10,2) circle (2pt) ; \fill (10,3) circle (2pt) ; \fill (10,5) circle (2pt) ; \fill (10,6) circle (2pt) ;
\fill (10,7) circle (2pt) ; \fill (10,8) circle (2pt) ; \fill (11,1) circle (2pt) ; \fill (11,2) circle (2pt) ; \fill (11,3) circle (2pt) ; \fill (11,5) circle (2pt) ; \fill (11,6) circle (2pt) ; \fill (11,7) circle (2pt) ; \fill (11,8) circle (2pt) ; \fill (12,1) circle (2pt) ; \fill (12,2) circle (2pt) ; \fill (12,3) circle (2pt) ; \fill (12,5) circle (2pt) ; \fill (12,6) circle (2pt) ; \fill (12,7) circle (2pt) ; \fill (12,8) circle (2pt) ; \fill (13,1) circle (2pt) ; \fill (13,2) circle (2pt) ; \fill (13,3) circle (2pt) ; \fill (13,5) circle (2pt) ; \fill (13,6) circle (2pt) ; \fill (13,7) circle (2pt) ; \fill (14,1) circle (2pt) ; \fill (14,2) circle (2pt) ; \fill (14,3) circle (2pt) ; \fill (14,5) circle (2pt) ; \fill (14,6) circle (2pt) ; \fill (14,7) circle (2pt) ; \fill (15,0) circle (2pt) ; \fill (15,1) circle (2pt) ; \fill (15,2) circle (2pt) ; \fill (15,3) circle (2pt) ; \fill (15,5) circle (2pt) ; \fill (15,6) circle (2pt) ; \fill (15,7) circle (2pt) ; \fill (16,4) circle (2pt) ; \fill (16,5) circle (2pt) ; \fill (16,6) circle (2pt) ; \fill (16,7) circle (2pt) ;
\draw (1,5) -- (1,8) -- (2,8) -- (2,5) 
      -- (3,5) -- (3,8) -- (4,8) -- (4,5)    
     --  (5,5) -- (5,7) -- (6,7) -- (6,5) -- (7,5) -- (7,7); \draw (16,4) -- (16,7) -- (15,7) -- (15,5) 
     -- (14,5) -- (14,7) -- (13,7)  -- (13,5)
     -- (12,5) -- (12,8) -- (11,8) -- (11,5)
     -- (10,5) -- (10,8) -- (9,8) -- (9,5)  ; \draw (7,0) -- (7,1) -- (6,1) -- (5,1) -- (4,1) 
      -- (4,2) -- (3,2) -- (3,1) -- (2,1) -- (2,2) -- (1,2) -- (1,1) ;     \draw (15,0) -- (15,3) -- (14,3) -- (14,1)  
     -- (13,1) -- (13,3) -- (12,3) -- (12,1)
     -- (11,1) -- (11,3) -- (10,3) -- (10,1) -- (9,1) -- (9,3) ;
\node at (0.7, 1) {\tiny 1} ;  \node at (0.7, 2) {\tiny 8} ;  \node at (0.7, 5) {\tiny 0} ;  \node at (0.7, 6) {\tiny 7} ;  \node at (0.7, 8) {\tiny 21} ;  \node at (4.3, 2) {\tiny 11} ;  \node at (4.3, 8) {\tiny 24} ;  \node at (7.3, 0) {\tiny 0} ;  \node at (7.3, 1) {\tiny 7} ;  \node at (7.3, 5) {\tiny 6} ;  \node at (7.3, 7) {\tiny 20} ;  \node at (8.7, 1) {\tiny 1} ;  \node at (8.7, 3) {\tiny 15} ;  \node at (8.7, 5) {\tiny 1} ;  \node at (8.7, 8) {\tiny 25} ;  \node at (12.3, 8) {\tiny 28} ; \node at (15.3, 0) {\tiny 0} ;  \node at (15.3, 1) {\tiny 7} ; \node at (15.3, 3) {\tiny 21} ; \node at (16.3, 4) {\tiny 0} ;  \node at (16.3, 5) {\tiny 8} ;  \node at (16.3, 7) {\tiny 24} ;  
\end{tikzpicture}
\end{center}
\end{figure}
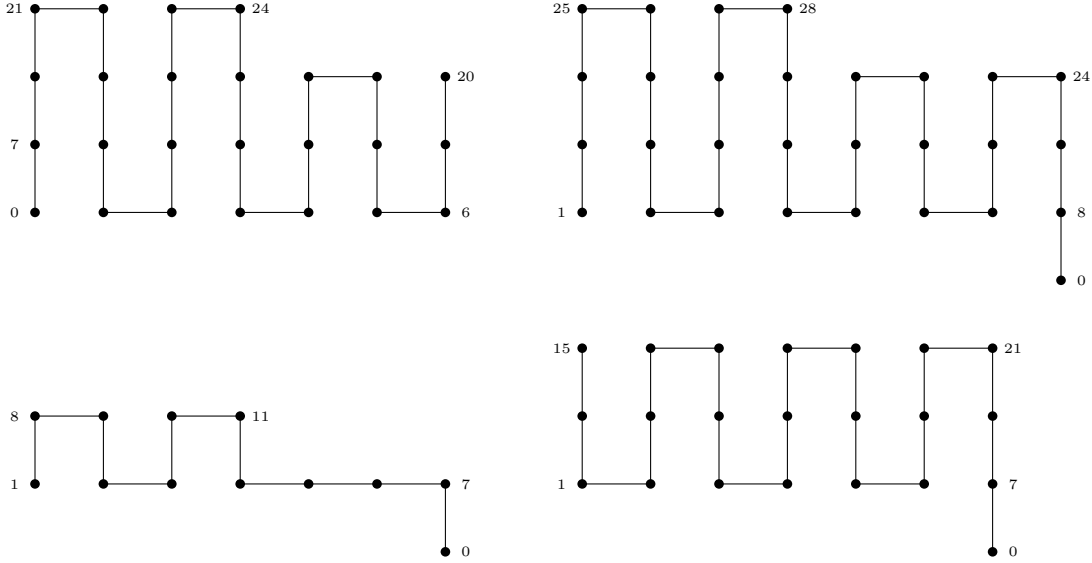

\begin{proof}[Proof Constructions.]
Trivial when~$b=0$. For~$b \geq 1$, we take the full fauxsets with respect to~$x$ and connect them with~$(x-1)$ 1-edges.  In particular, if~$r'$ is even, we take the fauxsets in increasing order from~0 (that is, $\varphi_0, \varphi_1, \ldots, \varphi_{x-1}$) and denote the construction~$\mathbf{h_1}$.  Otherwise, take the fauxsets in decreasing order from~0 (that is,  $\varphi_0, \varphi_{x-1}, \varphi_{x-2} \ldots, \varphi_{1}$) and denote the construction~$\mathbf{h_2}$.  Using the bridge concatenation notation, we have:
$$\mathbf{h_1} = \apph{}{i=0}{x-1} \Psi_i {\rm \ \ \ and \ \ \ }  \mathbf{h_2} =  \Psi_0 \uplus \left( \apph{}{i=x-1}{1} \Psi_i \right) . $$
Figure~\ref{fig:omega1} illustrates the four possibilities.
\end{proof}

\begin{lem}\label{lem:omega2}{\rm \cite{AO1} (Tail curl)}
There is a standard linear realization for~$\{1^x, x^b\}$.  
\end{lem}

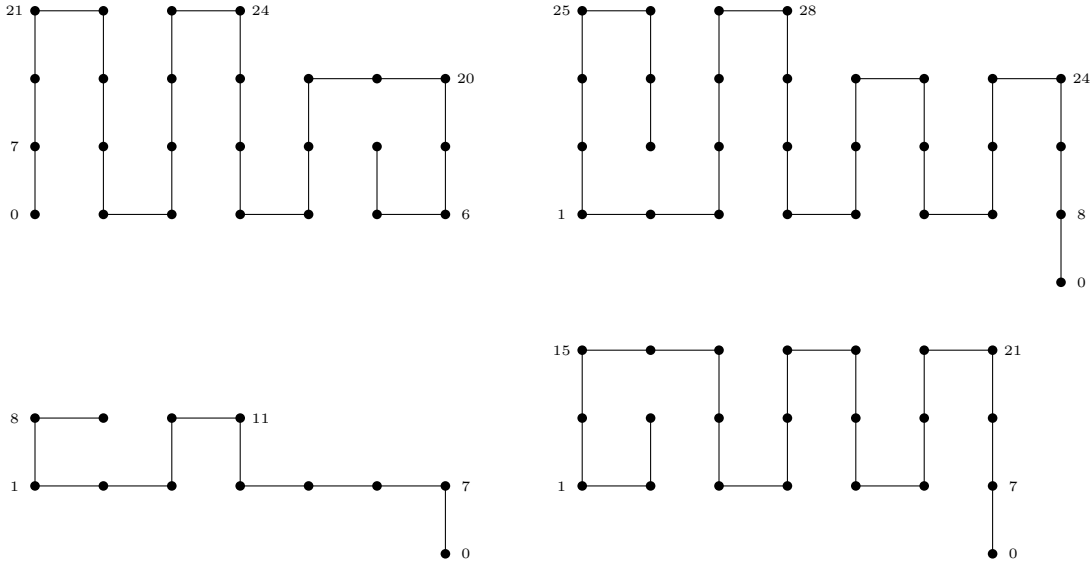
\begin{figure}
\caption{Standard linear realizations of~$\{1^{x}, x^b\}$ for $(v,x,b) = (25, 7,17)$, $(29, 8,20)$, $(12, 7, 4)$ and $(22, 7,14)$ obtained from those of Figure~\ref{fig:omega1} by applying a tail curl.}\label{fig:omega2}
\begin{center}
\begin{tikzpicture}[scale=0.9, every node/.style={transform shape}]
\fill (1,1) circle (2pt) ;\fill (1,2) circle (2pt) ;\fill (1,5) circle (2pt) ;\fill (1,6) circle (2pt) ;\fill (1,7) circle (2pt) ;\fill (1,8) circle (2pt) ;\fill (2,1) circle (2pt) ;\fill (2,2) circle (2pt) ;\fill (2,5) circle (2pt) ;\fill (2,6) circle (2pt) ;\fill (2,7) circle (2pt) ;\fill (2,8) circle (2pt) ;\fill (3,1) circle (2pt) ;\fill (3,2) circle (2pt) ;\fill (3,5) circle (2pt) ;\fill (3,6) circle (2pt) ;\fill (3,7) circle (2pt) ;\fill (3,8) circle (2pt) ;\fill (4,1) circle (2pt) ;\fill (4,2) circle (2pt) ;\fill (4,5) circle (2pt) ;\fill (4,6) circle (2pt) ;\fill (4,7) circle (2pt) ;\fill (4,8) circle (2pt) ;\fill (5,1) circle (2pt) ;\fill (5,5) circle (2pt) ;\fill (5,6) circle (2pt) ;\fill (5,7) circle (2pt) ;\fill (6,1) circle (2pt) ;\fill (6,5) circle (2pt) ;\fill (6,6) circle (2pt) ;\fill (6,7) circle (2pt) ;\fill (7,0) circle (2pt) ;\fill (7,1) circle (2pt) ;\fill (7,5) circle (2pt) ;\fill (7,6) circle (2pt) ;\fill (7,7) circle (2pt) ;\fill (9,1) circle (2pt) ;\fill (9,2) circle (2pt) ;\fill (9,3) circle (2pt) ;\fill (9,5) circle (2pt) ;\fill (9,6) circle (2pt) ;\fill (9,7) circle (2pt) ;\fill (9,8) circle (2pt) ;\fill (10,1) circle (2pt) ;\fill (10,2) circle (2pt) ;\fill (10,3) circle (2pt) ;\fill (10,5) circle (2pt) ;\fill (10,6) circle (2pt) ;\fill (10,7) circle (2pt) ;\fill (10,8) circle (2pt) ;\fill (11,1) circle (2pt) ;\fill (11,2) circle (2pt) ;\fill (11,3) circle (2pt) ;\fill (11,5) circle (2pt) ;\fill (11,6) circle (2pt) ;\fill (11,7) circle (2pt) ;\fill (11,8) circle (2pt) ;\fill (12,1) circle (2pt) ;\fill (12,2) circle (2pt) ;\fill (12,3) circle (2pt) ;\fill (12,5) circle (2pt) ;\fill (12,6) circle (2pt) ;\fill (12,7) circle (2pt) ;\fill (12,8) circle (2pt) ;\fill (13,1) circle (2pt) ;\fill (13,2) circle (2pt) ;\fill (13,3) circle (2pt) ;\fill (13,5) circle (2pt) ;\fill (13,6) circle (2pt) ;\fill (13,7) circle (2pt) ;\fill (14,1) circle (2pt) ;\fill (14,2) circle (2pt) ;\fill (14,3) circle (2pt) ;\fill (14,5) circle (2pt) ;\fill (14,6) circle (2pt) ;\fill (14,7) circle (2pt) ;\fill (15,0) circle (2pt) ;\fill (15,1) circle (2pt) ;\fill (15,2) circle (2pt) ;\fill (15,3) circle (2pt) ;\fill (15,5) circle (2pt) ;\fill (15,6) circle (2pt) ;\fill (15,7) circle (2pt) ;\fill (16,4) circle (2pt) ;\fill (16,5) circle (2pt) ;\fill (16,6) circle (2pt) ;\fill (16,7) circle (2pt) ;
\draw (1,5) -- (1,8) -- (2,8) -- (2,5)  -- (3,5) -- (3,8) -- (4,8) -- (4,5) --  (5,5) -- (5,7) -- (6,7) -- (7,7) -- (7,5) -- (6,5) -- (6,6); \draw (16,4) -- (16,7) -- (15,7) -- (15,5) -- (14,5) -- (14,7) -- (13,7)  -- (13,5)  -- (12,5) -- (12,8) -- (11,8) -- (11,5) -- (10,5) -- (9,5) -- (9,8) -- (10,8) -- (10,6) ; \draw (7,0) -- (7,1) -- (6,1) -- (5,1) -- (4,1)  -- (4,2) -- (3,2) -- (3,1) -- (2,1) -- (1,1) -- (1,2) -- (2,2) ;    \draw (15,0) -- (15,3) -- (14,3) -- (14,1)   -- (13,1) -- (13,3) -- (12,3) -- (12,1)   -- (11,1) -- (11,3) -- (10,3) -- (9,3) -- (9,1) -- (10,1) -- (10,2);
\node at (0.7, 1) {\tiny 1} ;  \node at (0.7, 2) {\tiny 8} ;  \node at (0.7, 5) {\tiny 0} ;  \node at (0.7, 6) {\tiny 7} ;  \node at (0.7, 8) {\tiny 21} ;  \node at (4.3, 2) {\tiny 11} ;  \node at (4.3, 8) {\tiny 24} ;  \node at (7.3, 0) {\tiny 0} ;  \node at (7.3, 1) {\tiny 7} ;  \node at (7.3, 5) {\tiny 6} ;  \node at (7.3, 7) {\tiny 20} ;  \node at (8.7, 1) {\tiny 1} ;  \node at (8.7, 3) {\tiny 15} ; \node at (8.7, 5) {\tiny 1} ;  \node at (8.7, 8) {\tiny 25} ;  \node at (12.3, 8) {\tiny 28} ; \node at (15.3, 0) {\tiny 0} ;  \node at (15.3, 1) {\tiny 7} ; \node at (15.3, 3) {\tiny 21} ;  \node at (16.3, 4) {\tiny 0} ;  \node at (16.3, 5) {\tiny 8} ;  \node at (16.3, 7) {\tiny 24} ;  
\end{tikzpicture}
\end{center}
\end{figure}

\begin{proof}[Proof Constructions.]
Trivial when~$b=0$. For~$b \geq 1$, consider the constructions~$\mathbf{h_1}$ and~$\mathbf{h_2}$ of Lemma~\ref{lem:omega1}. 
Let~$f$ be the final element of the linear realization.

Applying a {\em tail curl} refers to the following process: add the 1-edge from~$f$ to the previous fauxset to be traversed ($\varphi_{x-2}$ for~$\mathbf{h_1}$ and~$\varphi_2$ for~$\mathbf{h_2}$) and then, if possible, remove an $x$-edge internal to that fauxset in such a way that we have a Hamiltonian path.  
When possible, a tail curl transforms a linear realization for~$\{1^{x-1}, x^{b+1} \}$ to one for~$\{1^x , x^b\}$. 

If~$x$ and~$b$ meet one of the constraints of Lemma~\ref{lem:omega1}, we may obtain the desired linear realization for~$\{1^x, x^b\}$ via Lemma~\ref{lem:concat1} with~$s=1$.  Otherwise, we may take one of the linear realizations~$\mathbf{h_1}$ or~$\mathbf{h_2}$ for~$\{1^{x-1}, x^{b+1}\}$ and apply a tail curl.  (Often, both options are available.) 

Figure~\ref{fig:omega2} illustrates the application of a tail curl to each of the linear realizations of Figure~\ref{fig:omega1}.
\end{proof}

Observe how the diagrams in Figures~\ref{fig:omega1} and~\ref{fig:omega2} make the constructions of Lemmas~\ref{lem:omega1} and~\ref{lem:omega2} visually clear.  We draw them by first choosing one element~$x \in L$ with respect to which fauxsets were considered, and then aligning the vertices so that~$x$-edges (and hence fauxsets) are vertical and 1-edges are horizontal.  Following~\cite{AO1,AO2}, we call such representations {\em grid-based graphs}.  They are particularly useful for $\omega$-constructions, the grid-based graphs of which resemble battlements/crenellations with merlons and embrasures/crenels.

A common construction method of the paper can be thought of as taking an~$\omega$-construction and making alterations to the lowest row of~$x$-edges.  In such cases, we will draw our diagrams with only the lowest few rows.  To indicate how to reconstruct the unaltered remainder of the~$\omega$-construction from such a diagram, we will use and expand on existing terminology and also introduce some visual notation.

First, we formalize merlon extensions at the top of our constructions.  Suppose we are considering fauxsets with respect to~$y$ in~$K_v$. Say that a linear realization is of {\em type $\mathcal{C}_y$} in either of the following two cases:

\begin{itemize}
\item $y$ is even and there are edges between $v+1-i$ and~$v-i$ for all even~$i$ in the range~$2 \leq i \leq y$,
\item $y$ is odd, $v-1$ is an end-vertex, and there are edges between $v-i$ and~$v-i-1$ for all even~$i$ in the range~$2 \leq i \leq y-1$.
\end{itemize}

The first case corresponds exactly to the definition in~\cite{OPPS}; the second case is a new addition here.  In the language of~\cite{OPPS2}, a linear realization of type~$\mathcal{C}_y$ satisfies the definition of ``$y$-growable at $v-1$" while having additional properties.  We discuss $y$-growability further below.  Lemma~\ref{lem:cy} encapsulates how type $\mathcal{C}_y$ linear realizations can be used.

\begin{lem}\label{lem:cy} {\rm (Type $\mathcal{C}_y$)}
Suppose there is a type $\mathcal{C}_y$ linear realization for the multiset~$L$.   Then:

\begin{itemize}
\item there is a type $\mathcal{C}_y$ linear realization for the multiset~$L \cup \{ y^y \}$,
\item there is a linear realization for the multiset~$L \cup \{ y^{c'} \}$ for any even $c' < y$.
\end{itemize}

\end{lem}

\begin{proof}[Proof Construction.]
The first result is an application of~\cite[Theorem~2.4]{OPPS2}. When $y$ is even, both results follow from~\cite[Proposition~2.9]{OPPS}. When $y$ is odd, the second result is an easy generalization of it. 

Suppose~$L$ is a multiset in~$K_v$.  For the first item we need a linear realization in~$K_{v+y}$.  Obtain this by replacing each of the 1-edges between~$v-j$ and~$v-j+1$ guaranteed by the definition of type~$\mathcal{C}_y$ by the sequence 
$$(v-j, v+y-j, v+y-j+1, v-j+1),$$ 
adding two $y$-edges each time, and, if~$y$ is odd, add a $y$-edge from the final vertex~$v-1$ to~$v+y-1$.

The second item is similar.  Again, replace the guaranteed 1-edges with the same sequence.  Do so one at a time starting with the largest value of~$j$ and working down until the desired number of $y$-edges have been added. 
\end{proof}

Note that the constructions of Lemma~\ref{lem:cy} preserve whether the linear realization is standard, but the second item does not preserve whether it is perfect.

In the diagrams, we illustrate the appropriate 1-edges of a type~$\mathcal{C}_y$ by a thickened edge, with a small arrow pointing in the direction of increases by~$y$.  If~$y$ is odd, we  also add an arrow to the~$v-1$ vertex.

\begin{figure}
\caption{Examples of type $\mathcal{C}_y$ standard linear realizations for~$\{1^6, 7^{7}\}$ and~$\{1^8, 8^2\}$. }\label{fig:cy1}
\begin{center}
\begin{tikzpicture}[scale=0.9, every node/.style={transform shape}]
 \fill (15 ,0) circle (2pt);
\foreach \p in {0,...,6}{ \fill (\p ,1) circle (2pt);}
\foreach \p in {8,...,15}{ \fill (\p ,1) circle (2pt);}
\foreach \p in {0,...,6}{ \fill (\p ,2) circle (2pt);}
\foreach \p in {8,...,9}{ \fill (\p ,2) circle (2pt);}
\draw (0,1)--(0,2)--(1,2)--(1,1)--(2,1)--(2,2)--(3,2)--(3,1)--(4,1)--(4,2)--(5,2)--(5,1)--(6,1)--(6,2); 
\draw (15,0)--(15,1)--(8,1)--(8,2)--(9,2);
\foreach \p in {0,2} {\draw [very thick] (\p,2)--(\p+1,2) ; \draw [-stealth]  (\p+0.5,2) -- (\p+0.5, 2.3) ;}
\foreach \p in {4,8} {\draw [very thick] (\p,2)--(\p+1,2) ; \draw [-stealth]  (\p+0.5,2) -- (\p+0.5, 2.3) ;}
\foreach \p in {10,12,14} {\draw [very thick] (\p,1)--(\p+1,1) ; \draw [-stealth]  (\p+0.5,1) -- (\p+0.5, 1.3) ;}
\draw [-stealth]  (6,2) -- (6, 2.3);
\end{tikzpicture}
\end{center}
\end{figure}
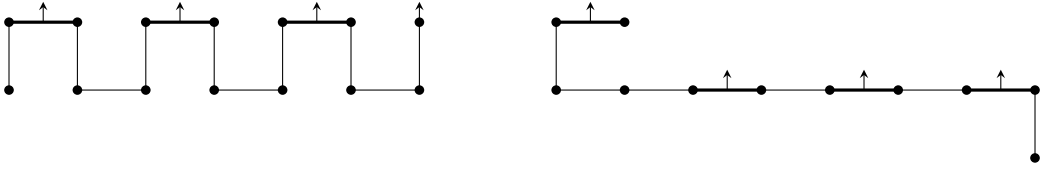

\begin{exa}\label{ex:C_y}{\rm (Type $\mathcal{C}_y$)} Write~$v = qy+r$ with $0 \leq r < y$.  

The first diagram of Figure~\ref{fig:cy1} illustrates the $\omega$-constructions for~$\{1^6, 7^b\}$ when~$b \geq 7$ and~$r$ is even. (An entirely horizontal line would also have worked and removed the condition on~$b$, but this would be less similar to our use cases for the paper.)  From this diagram, the first diagram of Figure~\ref{fig:omega1} can be reconstructed using one application of the first item of Lemma~\ref{lem:cy} and then one of the second to add~4 more 7-edges.

The second diagram of Figure~\ref{fig:cy1} illustrates the tail-curl constructions for~$\{1^8, 8^b\}$ when~$b \geq 2$ and~$r$ is odd.  (This time the diagram is flattened as much as possible.)  From this diagram, the second diagram of Figure~\ref{fig:omega2} can be reconstructed using two applications of the first item of Lemma~\ref{lem:cy} followed by an application of the second to add~4 more 8-edges.
\end{exa}

Next, we formalize lower merlon extensions.  These are the situations in which we need to insert $y$-edges internally to a diagram or at the bottom.  We do not have quite as much control in these cases and only accomplish something analogous to the first item in Lemma~\ref{lem:cy}.  

Call a linear realization~$\mathbf{h}$ that realizes~$L$ {\em $y$-growable at $m$} if we may construct a linear realization of $L \cup \{ y^y \}$ in $K_{v+y}$ via the following process.  Embed the subgraph on the labels $\{0,1,\ldots, m\}$ into~$K_{v+y}$ and also embed the subgraph on the labels $\{m-y+1, m-y+2, v-1\}$, translated up by~$y$ to use the labels $\{m+1, m+2, \ldots, v+y-1\}$, into~$K_{v+y}$.  Add $y$-edges between the labels~$m-y+i$ and~$m+i$ for each~$i$ in the range~$1\leq i \leq y$.  For each $x$-edge between two vertices in the set $\{m-y+1, m-y+2, \ldots, m\}$ in~$K_v$, there are now two $x$-edges in~$K_{v+y}$.  Remove one from each pair in such a way that a Hamiltonian path results.

This definition of $y$-growability is equivalent to that of~\cite{OPPS2} when applied to linear realizations.  The linear realization for $L \cup \{y^y\}$ obtained by {\em growing} a linear realization for~$L$ is also $y$-growable at~$m$~\cite{OPPS}.

Visually, in a grid-based diagram, we move the top portion of the picture one step in the $y$ direction, fill in the gaps with~$y$-edges, and remove duplicated edges in the obvious way.  We denote $y$-growability by circling the vertices $\{m-y+1, m-y+2, \ldots, m\}$ and indicating the $y$-direction with a small arrow on the vertex~$m$.  

A useful special case is when~$m = y-1$, and this can be thought of as quite like a type $\mathcal{C}_y$ situation, but at the bottom rather than the top.  An analogous approach to the second item of Lemma~\ref{lem:cy} could be considered in this case, but would not preserve the standardness of the linear realizations, so we do not pursue it here.

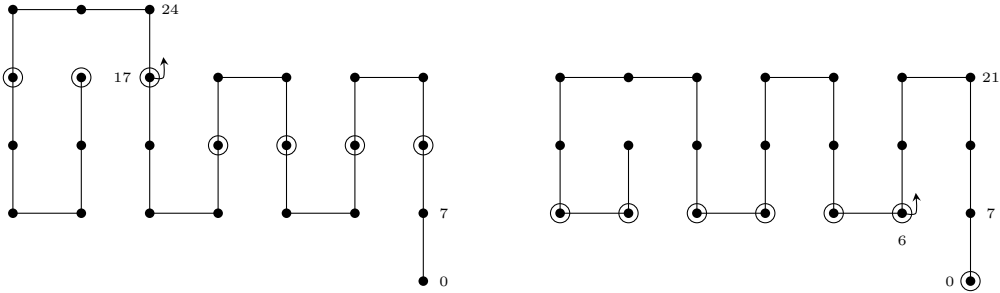
\begin{figure}
\caption{Examples of 7-growable standard linear realizations for $\{1^7, 7^{16}\}$ and $\{1^7, 7^{14}\}$. }\label{fig:grow}
\begin{center}
\begin{tikzpicture}[scale=0.9, every node/.style={transform shape}]
 \fill (6 ,0) circle (2pt);
  \fill (14 ,0) circle (2pt);
\foreach \p in {0,...,6}{ \fill (\p ,1) circle (2pt);}
\foreach \p in {8,...,14}{ \fill (\p ,1) circle (2pt);}
\foreach \p in {0,...,6}{ \fill (\p ,2) circle (2pt);}
\foreach \p in {8,...,14}{ \fill (\p ,2) circle (2pt);}
\foreach \p in {0,...,6}{ \fill (\p ,3) circle (2pt);}
\foreach \p in {8,...,14}{ \fill (\p ,3) circle (2pt);}
\foreach \p in {0,...,2}{ \fill (\p ,4) circle (2pt);}
\draw (6,0)--(6,3)--(5,3)--(5,1)--(4,1)--(4,3)--(3,3)--(3,1)--(2,1)--(2,4)--(0,4)--(0,1)--(1,1)--(1,3);
\draw (14,0)--(14,3)--(13,3)--(13,1)--(12,1)--(12,3)--(11,3)--(11,1)--(10,1)--(10,3)--(8,3)--(8,1)--(9,1)--(9,2);

\foreach \p in {0,1,2}{\draw (\p ,3) circle (4pt);}
\foreach \p in {3,...,6}{\draw (\p ,2) circle (4pt);}
\foreach \p in {8,...,13}{\draw (\p ,1) circle (4pt);}
\draw (14,0) circle (4pt);
\draw  [-stealth] plot [smooth] coordinates {(2,3) (2.2,3) (2.2,3.3)};
\draw  [-stealth] plot [smooth] coordinates {(13,1) (13.2,1) (13.2,1.3)};
\node at (6.3, 0) {\tiny 0} ;  \node at (6.3, 1) {\tiny 7} ;  \node at (2.3, 4) {\tiny 24} ;  \node at (1.6, 3) {\tiny 17} ;   \node at (13.7, 0) {\tiny 0} ;   \node at (14.3, 1) {\tiny 7} ;   \node at (14.3, 3) {\tiny 21} ;   \node at (13, 0.6) {\tiny 6} ;  
\end{tikzpicture}
\end{center}
\end{figure}

\begin{exa}\label{ex:y-growable}{\rm ($y$-Growability)} 
The first diagram in Figure~\ref{fig:grow} illustrates the 7-growability of the tail-curl $\omega$-construction for~$\{1^7, 7^{16}\}$ at~$m = 17$.  The second diagram  illustrates the 7-growability of the tail-curl construction for~$\{1^7, 7^{14}\}$ (given in the fourth diagram of Figure~\ref{fig:omega2}) at~$m = 6$.  From these we are able to construct the tail-curl constructions for~$\{1^7, 7^{7t+2}\}$ and $\{1^7, 7^{7t}\}$ for~$t > 2$.  Note that in each case there are many available choices of~$m$.
\end{exa}

\begin{figure}
\caption{A perfect linear realization for $\{ 5^7, 6^4\}$ and a standard linear realization for $\{ 1 , 4^4 , 6^6 \}$.    }\label{fig:slanted}
\begin{center}
\begin{tikzpicture}[scale=0.85, every node/.style={transform shape}]
\fill (0,3) circle (2pt) ; \fill (5.5,0) circle (2pt) ; 
\foreach \p in {1,...,5}{ \fill (\p ,1) circle (2pt);}
\foreach \p in {8,...,13}{ \fill (\p ,1) circle (2pt);}
\foreach \p in {0.5,1.5,...,4.5}{ \fill (\p ,2) circle (2pt);}
\foreach \p in {9,...,14}{ \fill (\p ,2) circle (2pt);}
\draw (0,3)--(1,1)--(1.5,2)--(2,1)--(2.5,2)--(3,1)--(3.5,2)--(4,1)--(4.5,2)--(5.5,0);
\draw (8,1)--(9,2)--(10,1)--(11,2)--(12,1)--(13,2)--(14,2)--(13,1)--(12,2)--(11,1)--(10,2)--(9,1);
\node at (-0.3, 3) {\tiny 11} ; \node at (0.7, 1) {\tiny 1} ; \node at (4.8, 2) {\tiny 10} ;  \node at (5.3, 1) {\tiny 5} ; \node at (5.8, 0) {\tiny 0} ; \node at (7.7, 1) {\tiny 0} ; \node at (8.7, 2) {\tiny 6} ; \node at (13.3, 1) {\tiny 5} ; \node at (14.3, 2) {\tiny 11} ;
\end{tikzpicture}
\end{center}
\end{figure}
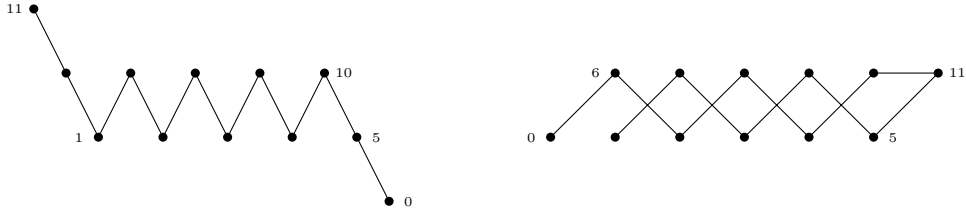

One last note on the visualizations: For many diagrams in this paper, we adjust the grid-based approach slightly by taking a vertical step of a realization with support $\{1, x, y\}$ to correspond to moving by $(x+y)/2$.  This slanted view gives the $x$ and $y$-edges equal visual prominence and more readily reveals the symmetry within and between some of our constructions.  When $(x+y)/2$ is not an integer, the vertices in different rows are offset. Figure~\ref{fig:slanted} illustrates this approach with supports~$\{5,6\}$ and~$\{1,4,6\}$.

\section{New useful standard linear realizations}\label{sec:para}

The overarching goal in this section is to construct useful standard linear realizations, which, when concatenated with altered $\omega$-constructions, yield larger linear realizations with as few 1-edges as possible.

We first build some perfect linear realizations for~$\{x^b, y^c\}$ when~$x$ and~$y$ are coprime.  When~$y=x+1$, this reduces to constructions of~\cite{AO1,OPPS}.  The first diagram of Figure~\ref{fig:slanted} gives an example of this special case when $\{x^b, y^c\} = \{5^7,6^4\}$. 

Second, we develop a method based on spanning linear forests and apply it to produce standard linear realizations for some multisets of the form~$\{1^{k-1}, (y-k)^b, y^c \}$ for small~$k$.

\begin{thm}\label{th:2perf}{\rm (Core perfect linear realizations)}
Let $x,y>1$ with $x$ and~$y$ coprime.  Then there is a perfect linear realization for $L = \{ x^{y-1}, y^{x+1} \}$.
\end{thm}

\begin{proof}
We have $v = (y-1) + (x+1) + 1 = x+y+1$.  We build the required realization using full fauxsets with respect to~$y$, with $x$-edges connecting them via bridge concatenations.  Write~$v = qy + r$ with $0 \leq r < y$.  

{\em Case~1:} $r = 0$.
Each fauxset has the same size, with fauxset~$\phi_i$ having smallest element~$i$ and largest element~$(q-1)y + i$.  Thus the length of the edge between the largest element in~$\phi_i$ and the smallest element in~$\phi_{i+1}$ is 
$$(q-1)y + i - (i+1) = qy - y  - 1 = v - y - 1 = x$$
and the construction
$$\apph{\mathit{x}}{k=0}{y-1} \Psi_{k} $$
is a Hamiltonian path that realizes~$L$.
This is illustrated for $(x,y) = (14,,5)$ in the first diagram of Figure~\ref{fig:perf}.

\begin{figure}
\caption{Perfect linear realizations for~$\{5^{15}, 14^4\}$ and $\{7^{10},9^6\}$.  }\label{fig:perf}
\begin{center}
\begin{tikzpicture}[scale=0.9, every node/.style={transform shape}]
\foreach \p in {0,...,4}{ \fill (\p ,0) circle (2pt);}
\foreach \p in {0,...,4}{ \fill (\p ,1) circle (2pt);}
\foreach \p in {0,...,4}{ \fill (\p ,2) circle (2pt);}
\foreach \p in {0,...,4}{ \fill (\p ,3) circle (2pt);}
\foreach \p in {6,...,12}{ \fill (\p ,1) circle (2pt);}
\foreach \p in {6,...,12}{ \fill (\p ,2) circle (2pt);}
\foreach \p in {6,...,8}{ \fill (\p ,3) circle (2pt);}
\draw (0,0)--(0,3)--(1,0)--(1,3)--(2,0)--(2,3)--(3,0)--(3,3)--(4,0)--(4,3);
\draw (6,1)--(6,3);
\draw (7,3)--(7,1)--(9,2)--(9,1)--(11,2)--(11,1);
\draw (8,3)--(8,1)--(10,2)--(10,1)--(12,2)--(12,1);
\draw  plot [smooth] coordinates {(6,3) (5.5,0) (11,1)};
\draw  plot [smooth] coordinates {(7,3) (6.5,0) (12,1)};
\node at (-0.3, 0) {\tiny 0} ;   
\node at (-0.3, 1) {\tiny 5} ;   
\node at (-0.3, 3) {\tiny 15} ;   
\node at (4.3, 0) {\tiny 4} ;  
\node at (4.3, 3) {\tiny 19} ;  
\node at (6, 0.7) {\tiny 0} ;   
\node at (6.3, 2) {\tiny 7} ;   
\node at (6.3, 3) {\tiny 14} ;   
\node at (8.3, 3) {\tiny 16} ;   
\node at (12.3, 1) {\tiny 6} ;   
\node at (12.3, 2) {\tiny 13} ;   
\end{tikzpicture}
\end{center}
\end{figure}
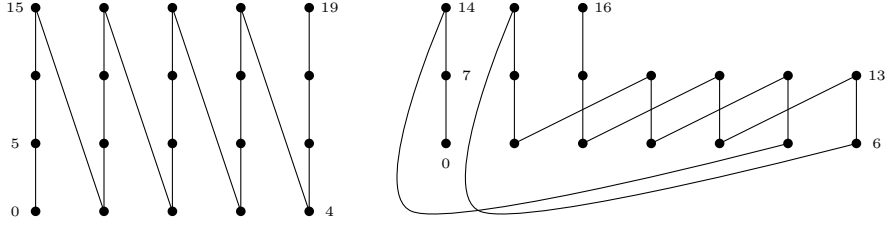

{\em Case~2:} $r > 0$.
We have that $v \equiv x+1 \pmod{y}$ and the largest element in fauxset~$\phi_i$ is $qy + i$ for~$i < r$ and~$(q-1)y + i$ otherwise.  Let~$t = qy-x = y-r+1$ and note that~$\gcd(t,y) = 1$.  Let $[g_0, g_1, \ldots, g_{y-1}]$ be a sequence $[0, t, 2t, \ldots, (y-1)t]$ in~$\Z_y$.  As $\gcd(t,y) = 1$, this includes each element of $\Z_y$ exactly once.  Consider the construction
$$\apph{\mathit{x}}{i=0}{y-1} \Psi_{g_i},$$
starting at the vertex~0.
This has the required differences by construction, so we only need to check that it is a valid Hamiltonian path.  

As in the previous case, all bridging edges are between the largest vertex in one fauxset and the smallest one in another (it may be that there are fauxsets with a single element that plays both roles).  There are two sizes of fauxset; they have either~$q+1$ or~$q$ elements.

The bridge between $\Psi_{g_i}$ and~$\Psi_{g_{i+1}}$ when~$\phi_{g_i}$ is one of the larger fauxsets is 
$$(qy+g_i) - g_{i+1} = (qy + g_i) - (g_i + t) = qy - (qy-x) = x.$$  
When~$\phi_{g_i}$ is one of the smaller fauxsets we have 
$$((q-1)y + g_i) - g_{i+1} = (qy - y + g_i) - (g_i + t - y) = x.$$
Therefore, the bridges have the required lengths and the construction is valid.

This is illustrated for $(x,y) = (9,7)$ in the second diagram of Figure~\ref{fig:perf}.
\end{proof}

Unusually for results of this type, Theorem~\ref{th:2perf} does not require that~$x<y$.  Corollary~\ref{cor:2perf} translates it for this situation, emphasizing that for a given pair of coprime numbers there are two perfect linear realizations available from the construction.

\begin{cor}\label{cor:2perf}{\rm (Split-form)}
Let~$1 < x < y$. Then there are perfect linear realizations for the multisets  $\{ x^{y-1}, y^{x+1} \}$ and $\{ x^{y+1}, y^{x-1} \}$.
\end{cor}

Here is another view on the relationship between the two perfect linear realizations of Corollary~\ref{cor:2perf}.  Suppose the perfect linear realization for $\{ x^{y-1}, y^{x+1} \}$ is $[ h_1, h_2, \ldots, h_v ]$.  Observe that we have~$h_1 = 0$, $h_2 = y$, $h_{v-1} = v-1-y$ and $h_v = v-1$.   Rewriting in terms of~$x$, using that $v = x+y+1$, we obtain $h_2 = v-1-x$ and $h_{v-1} = x$.  Therefore the Hamiltonian path $[h_1, h_{v-2}, h_{v-3}, \ldots, h_2, h_{v-1}]$ obtained by reversing the order of the interior vertices, realizes $( \{ x^{y-1}, y^{x+1} \} \setminus \{ y^2 \}) \cup \{x^2\}  = \{ x^{y+1}, y^{x-1} \}$.  This is the other perfect linear realization from Corollary~\ref{cor:2perf}.

Whenever we have a perfect linear realization we may obtain a standard linear realization for a multiset with one fewer element by removing the final edge.  This yields Corollary~\ref{cor:2std}.

\begin{cor}\label{cor:2std}{\rm (Final edge removed)}
Let~$1 < x < y$. There are standard linear realizations for the multisets  $\{ x^{y-1}, y^{x} \}$ and $\{ x^{y}, y^{x-1} \}$.
\end{cor}

We now move to the second major construction of the section.
As with the $\omega$-constructions, the idea is to construct a spanning linear forest of~$K_v$ using relatively few paths and add 1-edges to make a Hamiltonian path.   This time we use both $x$-edges and~$y$-edges in the paths of the forest.

The spanning linear forest has~$k$ paths.  We construct it in two stages: first we define a spanning linear forest with~$k$ paths in $K_{2y}$ that linearly realizes~$\{(y-k)^{y-k} , y^y\}$ and then show that these can be combined to create a spanning linear forest linearly realizing $\{(y-k)^{j(y-k)} , y^{j y}\}$ for any~$j > 1$.

Let~$\Gamma^{(1)}_{k,y}$ be the subgraph of~$K_{2y}$ that includes the $y$-edges 
$$(0,y), (1,y+1),  \ldots (y-1, 2y-1)$$
(in fact, these are simply the fauxsets with respect to~$y$)
and the $(y-k)$-edges
$$(k, y), (k+1, y+1), \ldots, (y-1, 2y-k-1).$$  
Together, these form~$k$ paths, each with one end-point in $\{0, 1, \ldots, k-1\}$ and one in $\{2y-k, 2y-k+1, \ldots,  2y-1\}$.  So $\Gamma^{(1)}_{k,y}$ is a spanning linear forest with~$k$ paths that linearly realizes~$\{(y-k)^{y-k} , y^y\}$.

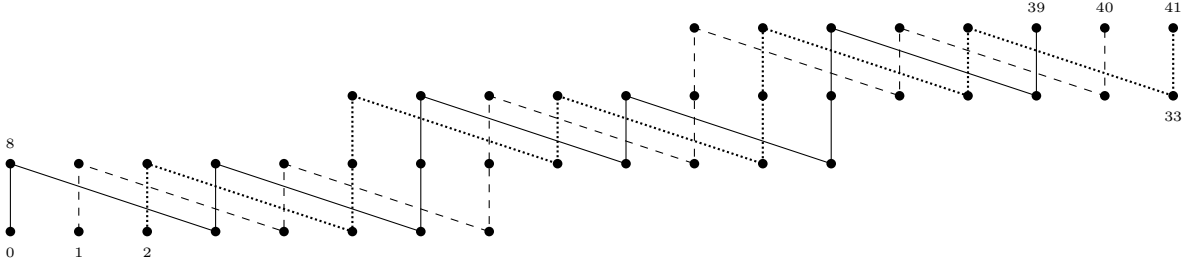
\begin{figure}
\caption{The spanning linear forest~$\Gamma^{(3)}_{3,8}$ that linearly realizes~$\{5^{15}, 8^{24}\}$ in $K_{42}$.  Different line patterns are used for the different paths to make the construction clearer.  Adding 1-edges between~39 and~40 and between~1 and~2 gives a perfect linear realization for~$\{1^2,5^{15}, 8^{24}\}$. }\label{fig:tracks}
\begin{center}
\begin{tikzpicture}[scale=0.9, every node/.style={transform shape}]
\foreach \p in {0,1,2,3,4,5,6,7}{\fill (\p ,0) circle (2pt) ;}
\foreach \p in {0,1,2,3,4,5,6,7,8,9,10,11,12}{\fill (\p ,1) circle (2pt) ;}
\foreach \p in {5,6,7,8,9,10,11,12,13,14,15,16,17}{\fill (\p ,2) circle (2pt) ;}
\foreach \p in {10,11,12,13,14,15,16,17}{\fill (\p ,3) circle (2pt) ;}
\draw (0,0) -- (0,1) -- (3,0) -- (3,1) -- (6,0) -- (6,2) -- (9,1) -- (9,2) -- (12,1) -- (12,3) -- (15,2) -- (15,3)  ;
\draw[dashed] (1,0) -- (1,1) -- (4,0) -- (4,1) -- (7,0) -- (7,2) -- (10,1) -- (10,3) -- (13,2) -- (13,3) -- (16,2) -- (16,3)  ;
\draw[thick, densely dotted] (2,0) -- (2,1) -- (5,0) -- (5,2) -- (8,1) -- (8,2) -- (11,1) -- (11,3) -- (14,2) -- (14,3) -- (17,2) -- (17,3)  ;
\node at (0, -0.3) {\tiny 0} ;  \node at (1, -0.3) {\tiny 1} ;  \node at (2, -0.3) {\tiny 2} ;  \node at (15, 3.3) {\tiny 39} ;  \node at (16, 3.3) {\tiny 40} ;  \node at (17, 3.3) {\tiny 41} ;  \node at (0, 1.3) {\tiny 8} ;  \node at (17, 1.7) {\tiny 33} ; 
\end{tikzpicture}
\end{center}
\end{figure}

In general, $\Gamma^{(j)}_{k,y}$ will be a spanning linear forest that linearly realizes $\{(y-k)^{j(y-k)} , y^{j y}\}$.  We describe how to construct $\Gamma^{(i+1)}_{k,y}$ from $\Gamma^{(i)}_{k,y}$.  Assume that $\Gamma^{(i)}_{k,y}$ is a spanning linear forest for~$K_v$ with~$k$ points having one end-point in $\{0, 1, \ldots, k-1\}$ and one in  $\{v-k, v-k+1, \ldots , v-1\}$.  The process is analogous to concatenation.  Take the translation of $\Gamma^{(1)}_{y,k}$ by $v-k$ and identify the end-point vertices $\{v-k, v-k+1, \ldots , v-1\}$ of $\Gamma^{(i)}_{k,y}$ and the translation of $\Gamma^{(1)}_{y,k}$.  It is straightforward to check that $\Gamma^{(i+1)}_{k,y}$ has the required properties and that it allows the continuation of the inductive definition.
Figure~\ref{fig:tracks} illustrates $\Gamma^{(3)}_{3,8}$.

Observe that the upper end-points of the paths occur in the same order as the lower end-points, up to a cyclic permutation.  That is, if we label the path that starts with point~$i$, for $0 \leq i < k$, as~$\mathbf{p_i}$, then the respective end-points $(v-k, v-k+1, \ldots , v-1)$ appear in that order, possibly with a cyclic shift.

This means that when $k \mid y$ the paths of~$\Gamma^{(j)}_{k,y}$ line up perfectly to give a standard linear realization, which gives the following result.

\begin{thm}\label{th:tracks_div}{\rm (When $k \mid y$)}
Suppose~$y = mk$ for some~$m > 1$.  For all~$j \geq 1$, there is a standard linear realization for $L = \{ 1^{k-1} , (y-k)^{j(y-k)} , y^{j y} \}$.  When~$k$ is odd there is a perfect linear realization.  
\end{thm}

\begin{proof}
The pairs of end-points of the~$k$ paths in $\Gamma^{(j)}_{k,mk}$ have the form $(i, v-k+i)$ for $0 \leq i < k$.  Labeling the path with end-point~$i$, for $0 \leq i < k$, as~$\mathbf{p_i}$, the bridge concatenation 
$$\apph{}{i=0}{k-1}  \mathbf{p_i}$$
starting at~0 has the required properties.
\end{proof}

The second diagram of Figure~\ref{fig:slanted} illustrates Theorem~\ref{th:tracks_div} for~$k=2$, $m=3$, and~$j=1$.

For very small~$k$, we do not need divisibility to obtain standard linear realizations.

\begin{thm}\label{th:tracks23}{\rm (Small $k$)}
Let~$k \in \{2,3,4\}$ and $y > k+1$.  For all~$j \geq 1$, there is a standard linear realization for $L = \{ 1^{k-1} , (y-k)^{j(y-k)} , y^{j y} \}$.  
\end{thm}

\begin{proof}
If the ordering of the upper end-points is $(v-k, v-k+1, \ldots, v-1)$ then we can use the method of Theorem~\ref{th:tracks_div}.

Otherwise, for~$k=2$ we have that the ordering of the  upper endpoints is~$(v-1, v-2)$.   Hence, using the same path notation as the proof of Theorem~\ref{th:tracks_div}, the bridge concatenation $\mathbf{p_0} \uplus \mathbf{p_1}$ is a valid construction and gives the required realization.

For~$k=3$ we have two possibilities for the order of the upper end-points: $(v-2, v-1, v-3)$ and $(v-1, v-3, v-2)$.  In both cases,  $\mathbf{p_0} \uplus \mathbf{p_2} \uplus \mathbf{p_1}$ gives the required realization; in the latter case it is perfect.

For~$k=4$ we have three possibilities for the order of the upper end-points: $(v-3, v-2, v-1, v-4)$, $(v-2, v-1, v-4, v-3)$ and $(v-1, v-4, v-3, v-2)$.  In all cases,  $\mathbf{p_0}  \uplus \mathbf{p_3} \uplus \mathbf{p_2} \uplus \mathbf{p_1}$ gives the required realization.
\end{proof}

The method of proof of Theorem~\ref{th:tracks23} breaks down at~$k=5$.  For example, the upper end-points might be ordered as $(v-2, v-1, v-5, v-4, v-3)$ for which there is no successful bridge concatenation using 1-edges.

\section{Revisiting the case of $k=1$}\label{sec:1y-1y}

The purpose of this section is twofold.  First, we recall the main techniques of~\cite{AO1} used in the construction of standard linear realizations for the case of $k=1$, as they will be used and generalized in the remaining sections of the paper.  Second, we introduce new constructions when~$b$ and~$c$ are small, with~$a \geq y - \min(b,c)$ (sometimes even smaller), covering infinitely many new families for arbitrarily large~$y$, and showing that the remaining possible exceptions to the Coprime BHR Conjecture when~$y \leq 16$ are realizable, proving Theorem~\ref{th:small_y-1}.

The constructions for~$k=1$ that were introduced in~\cite{AO1} use {\em $\gamma$-moves}, each moving a single vertex elsewhere on the same path.  The two particular $\gamma$-moves used in these constructions have the benefit of changing no more than~2 of the elements of the realized multiset, by replacing either~$y$-edges with~($y-1$)-edges, or vice versa. 

The slanted grid-based representation (with a vertical step corresponding to a distance of~$y-\frac{1}{2}$) illustrates the geometry of these two~$\gamma$-moves that are vertical reflections of one another:

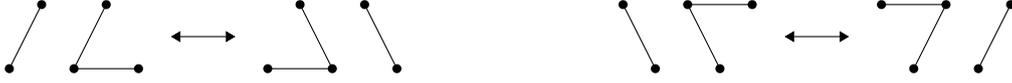
\begin{center}
\begin{tikzpicture}[scale=0.85, every node/.style={transform shape}]
\foreach \p in {0,1,2}{ \fill (\p ,0) circle (2pt);}
\foreach \p in {4,5,6}{ \fill (\p ,0) circle (2pt);}
\foreach \p in {0.5, 1.5}{ \fill (\p ,1) circle (2pt);}
\foreach \p in {4.5, 5.5}{ \fill (\p ,1) circle (2pt);}
\draw (0,0)--(0.5,1);
\draw (1.5,1)--(1,0)--(2,0);
\draw (4.5,1)--(5,0)--(4,0);
\draw (6,0)--(5.5,1);
\draw [Triangle-Triangle] (2.5,0.5) -- (3.5,0.5) ;
\foreach \p in {10,11}{ \fill (\p ,0) circle (2pt);}
\foreach \p in {14,15}{ \fill (\p ,0) circle (2pt);}
\foreach \p in {9.5, 10.5, 11.5}{ \fill (\p ,1) circle (2pt);}
\foreach \p in {13.5, 14.5, 15.5}{ \fill (\p ,1) circle (2pt);}
\draw (9.5,1)--(10,0);
\draw (11.5,1)--(10.5,1)--(11,0);
\draw (13.5,1)--(14.5,1)--(14,0);
\draw (15.5,1)--(15,0);
\draw [Triangle-Triangle] (12,0.5) -- (13,0.5) ;
\end{tikzpicture}
\end{center}

In order to replace only a single edge, two other types of moves were introduced in~\cite{AO1}: {\em corner cut} and {\em corner flap}, both involving the vertex in the path with the largest label.

Using $\gamma$-moves and the corner cut to replace~$y$-edges with ($y-1$)-edges in $\omega$-constructions for support~$\{1,y\}$, we get standard linear realizations when~$a \geq y$, for small~$b$:

\begin{lem}\label{lem:1_small_b}{\rm \cite{AO1} ($k=1$; small $b$)} 
Let $L = \{ 1^a, (y-1)^b, y^c \}$ with~$b < y-1$.  There is a standard linear realization for~$L$ when~$a \geq y$.
\end{lem}

Using $\gamma$-moves and the corner flap to replace ($y-1$)-edges with~$y$-edges in $\omega$-constructions for support~$\{1,y-1\}$, we get standard linear realizations when~$a \geq y$, for small~$c$:

\begin{lem}\label{lem:1_small_c}{\rm \cite{AO1} ($k=1$; small $c$)} 
Let $L = \{ 1^a, (y-1)^b, y^c \}$ where~$c < y$.  There is a standard linear realization for~$L$ when~$a \geq y$.
\end{lem}

Building a perfect linear realization for $\{ (y-1)^{j(y-1)}, y^{jy} \}$ by concatenating~$j$ copies of the perfect linear realization  $\{ (y-1)^{(y-1)}, y^{y} \}$ given by Theorem~\ref{th:2perf}, and concatenating it with the standard linear realization from the appropriate choice from  Lemma~\ref{lem:1_small_b} or~\ref{lem:1_small_c}, we get a standard linear realization when~$a \geq y$, for any multiset in the $k=1$ case:

\begin{thm}\label{th:y-1}{\rm \cite{AO1} ($k=1$; when $a \geq y$)}
Let $L = \{ 1^a, (y-1)^b, y^c \}$.  There is a standard linear realization for~$L$ whenever~$a \geq y$. 
\end{thm}

Using modular arithmetic on the multiplicities of equivalent multisets yields that the Coprime BHR~Conjecture holds in the case of~$k=1$ when~$v$ is large enough:

\begin{thm}\label{th:large_v}{\rm \cite{AO1} ($k=1$; large $v$)} 
If $v \geq 2y^2 + 9y$, then the Coprime BHR~Conjecture holds for multisets with support~$\{1,y-1,y\}$. 
\end{thm}

Applying the above results in conjunction with others from the existing literature yields a comprehensive result for the Coprime BHR Conjecture when~$y$ is small enough:

\begin{cor}\label{cor:small_y}{\rm \cite{AO1} ($k=1$; small $y$)}
If $y \leq 16$, then the Coprime BHR~Conjecture holds for multisets with support~$\{1,y-1,y\}$, except possibly for the following values of~$v$ and~$y$:

$$
\begin{array}{r||c|c|c|c|c|c}
y & 9 & 12 & 13 & 14 & 15 & 16  \\
\hline 
v & 43 & 41,53,79 & 41,67,89 & 67 & 73 & 41, 43, 73, 103, 137, 167  \\
\end{array}
$$

\end{cor}

We now introduce new constructions for when both~$b$ and~$c$ are small, requiring fewer 1-edges than known constructions.  These constructions dovetail with the existing methods, allowing us to show that the Coprime BHR Conjecture holds for~$k=1$ when~$y \leq 16$, thus proving Theorem~\ref{th:small_y-1}.

\begin{lem}\label{lem:sawtooth}{\rm ($k=1$; small $b$ and $c$)}
Let~$L = \{1^a, (y-1)^b, y^c\}$ with $0 < b,c < y-1$.  There is a standard linear realization for~$L$ in each of the following cases:
\begin{itemize}
\item $c > b-1$ and $a \geq y-b-1$,
\item $c = b-1$ and $a \geq y-c $,
\item $c < b-1$ and $a \geq y-c-2$.
\end{itemize}
In particular, there is a a standard linear realization for~$L$ when~$a \geq y - \min(b,c)$.
\end{lem}

\begin{proof}
The rough overarching idea is to take whichever of~$b$ and~$c$ is larger and think in terms of fauxsets with respect to~$y-1$ or~$y$ respectively, using the other value to transition between some of the fauxsets.

{\em Case 1:} $c > b-1$, $a = y-b-1$. If $c-b$ is even,  we take fauxsets with respect to~$y$ and use the standard linear realization 
$$\apph{\boldell}{k=0}{y-1}  \Psi_k \ \ \  {\textrm where } \ \ \ 
\boldell = (\underbrace{1,1, \ldots, 1}_{{c-b \ {\textrm terms}}},
\underbrace{ y-1,y-1, \ldots, y-1 }_{b \ {\textrm terms} }, 
\underbrace{ 1,1, \ldots, 1}_{y-c-1 \ {\textrm terms}  }  ).   $$
As~$c-b$ is even, each fauxset that we shall exit via a $(y-1)$-edge is being entered at the lowest value.

The last vertex of the above realization is~$y-1$.  We may add one more~$y$-edge, thus covering the case with~$c-b$ odd, by adding a $y$-edge from~$y-1$ to~$-1$ and then adding one to every vertex label.  

The first diagram of Figure~\ref{fig:sawtooth} illustrates this for~$\{1^5,7^2,8^6 \}$ and~$\{1^5,7^2,8^7 \}$.

\begin{figure}
\caption{Standard linear realizations for $\{1^5,7^2,8^6 \}$, $\{1^6,8^4,9^3 \}$, $\{1^3,8^7,9^6 \}$ and $\{1^5,8^6,9^2 \}$ (ignoring the dashed lines), and also $\{1^5,7^2,8^7 \}$ and $\{1^5,8^7,9^2 \}$ (including the dashed lines and adding~1 to the vertex labels) from the proof of Lemma~\ref{lem:sawtooth}.}\label{fig:sawtooth}
\begin{center}
\begin{tikzpicture}[scale=0.85, every node/.style={transform shape}]
\foreach \p in {1,...,8}{ \fill (\p ,1) circle (2pt);}
\foreach \p in {10,...,17}{ \fill (\p ,1) circle (2pt);}
\foreach \p in {0.5,1.5,...,7.5}{ \fill (\p ,2) circle (2pt);}
\foreach \p in {9.5,10.5,...,14.5}{ \fill (\p ,2) circle (2pt);}
\foreach \p in {1,...,8}{ \fill (\p ,4.5) circle (2pt);}
\foreach \p in {10,...,17}{ \fill (\p ,4.5) circle (2pt);}
\foreach \p in {1.5,2.5,...,6.5}{ \fill (\p ,5.5) circle (2pt);}
\foreach \p in {9.5,10.5,...,13.5}{ \fill (\p ,5.5) circle (2pt);}
\fill (8.5,0) circle (2pt); \fill (17.5,0) circle (2pt); \fill (7.5,3.5) circle (2pt); \fill (17.5,3.5) circle (2pt);
\draw (1,4.5)--(1.5,5.5)--(2.5,5.5)--(2,4.5)--(3,4.5)--(3.5,5.5)--(4.5,5.5)--(4,4.5)--(5,4.5)--(5.5,5.5)--(6,4.5)--(6.5,5.5)--(7,4.5)--(8,4.5);
\draw [dashed] (8,4.5)--(7.5,3.5);
\draw (10,4.5)--(9.5,5.5)--(11.5,5.5)--(11,4.5)--(12,4.5)--(12.5,5.5)--(13,4.5)--(13.5,5.5)--(14,4.5)--(17,4.5)--(17.5,3.5);
\draw (1,1)--(0.5,2)--(2.5,2)--(2,1)--(3,1)--(3.5,2)--(4,1)--(4.5,2)--(5,1)--(5.5,2)--(6,1)--(6.5,2)--(7,1)--(7.5,2)--(8.5,0);
\draw (10,1)--(9.5,2)--(10.5,2)--(11,1)--(12,1)--(11.5,2)--(12.5,2)--(13,1)--(13.5,2)--(14,1)--(14.5,2)--(15,1)--(17,1);
\draw [dashed] (17,1)--(17.5,0);
\node at (0.7, 1) {\tiny 1} ;  \node at (0.2, 2) {\tiny 9} ;  \node at (0.7, 4.5) {\tiny 0} ; \node at (1.2, 5.5) {\tiny 8} ;\node at (8.7, 0) {\tiny 0} ;    \node at (8.3, 1) {\tiny 8} ;  \node at (7.8, 2) {\tiny 16} ;  \node at (7.8, 3.5) {\tiny -1} ; \node at (8.3, 4.5) {\tiny 7} ;  \node at (6.8, 5.5) {\tiny 14} ;  \node at (9.7, 1) {\tiny 0} ;  \node at (9.2, 2) {\tiny 8} ;  \node at (9.7, 4.5) {\tiny 1} ; \node at (9.2, 5.5) {\tiny 9} ; \node at (17.8, 0) {\tiny -1} ; \node at (17.3, 1) {\tiny 7} ;  \node at (14.8, 2) {\tiny 14} ;  \node at (17.8, 3.5) {\tiny 0} ;   \node at (17.3, 4.5) {\tiny 8} ;  \node at (13.8, 5.5) {\tiny 14} ; 
\end{tikzpicture}
\end{center}
\end{figure}
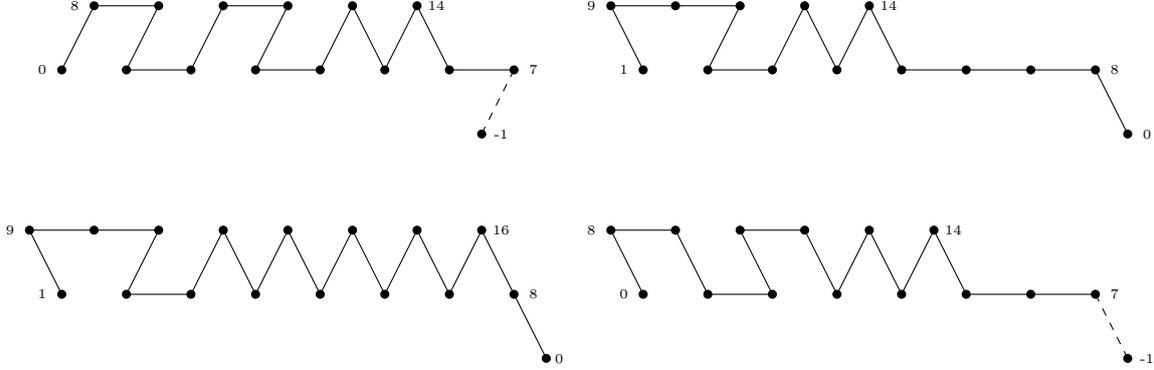

{\em Case 2:} $c = b-1$, $a = y-c$.  We take fauxsets with respect to~$y-1$.  The last~6 vertices will always be the sequence $\Theta = (3,2,y+2,y+1,y,1)$, which realizes~$\{1^3,y-1,y\}$ while traversing all elements of the fauxsets~$\Psi_1$, $\Psi_2$ and~$\Psi_3$.  
The standard linear realization is
$$\Psi_0 \ \uplus^{i} \ \left(\apph{\boldell}{k=y-2}{4}  \Psi_k \right) \uplus \ \Theta$$
where~$i=y$ when~$b=y-2$ and~$i=1$ otherwise and 
$$\boldell = (\underbrace{ 1, 1, \ldots, 1}_{y-c-3 \ {\textrm terms}} , 
    \underbrace{ y, y,  \ldots, y}_{c-1 \ {\textrm terms}})$$
Similarly to the previous case, each fauxset that we shall exit via either a $y$-edge or the sequence~$\Theta$ is being entered at the lowest element.  

The second and third diagrams of Figure~\ref{fig:sawtooth} illustrate this for~$\{1^6,8^4,9^3\}$ and~$\{1^3,8^3,9^6\}$, which use~~$i=1$ and~$i=y=9$, respectively.

{\em Case 3:} $c < b-1$, $a = y- c- 2$. If~$b-c$ even, take fauxsets with respect to~$y-1$ and use the standard linear realization
$$\apph{\boldell}{k=0}{y-2}  \Psi_k \ \ \  {\textrm where } \ \ \ 
\boldell = (\underbrace{1,1, \ldots, 1}_{{b-c-1 \ {\textrm terms}}},
\underbrace{ y-1,y-1, \ldots, y-1 }_{c \ {\textrm terms} }, 
\underbrace{ 1,1, \ldots, 1}_{y-b-1 \ {\textrm terms}  }  ).   $$
As~$b-c-1$ is odd, we enter the first (and each subsequent) fauxset from which we are going to transition with a $y$-edge at the highest value.

The last vertex of the above realization is~$y-2$.  We may add one more~$(y-1)$-edge, thus covering the case with~$b-c$ odd, by adding a $(y-1)$-edge from~$y-2$ to~$-1$ and then adding one to every vertex label.  

The last diagram of Figure~\ref{fig:sawtooth} illustrates this for~$\{1^5,8^5,9^2\}$ and~$\{1^5,8^6,9^2\}$.

As all the linear realizations constructed are standard, we may prepend a perfect linear realization with support~$\{1\}$ of the appropriate length to obtain standard linear realizations with any larger number of 1-edges.

The final statement follows immediately, with possible equality from the case~$c=b-1$.
\end{proof}

Using Lemma~\ref{lem:sawtooth} we are able resolve the possible exceptions to the Coprime BHR Conjecture noted in Corollary~\ref{cor:small_y}.

\begin{thm}\label{th:15cases}{\rm ($k=1$; small $y$)} 
For each outstanding case for $y \leq 16$, the Coprime BHR~Conjecture holds for multisets with support~$\{1,y-1,y\}$. 
\end{thm}

\begin{proof}
We work through each potential counterexample and show that all are realizable.

{\em Case~1:} $y=9$, $v=43$.  The three equivalent multisets are~$L = \{1^a,8^b,9^c\}$, $L' = \{ 1^b,15^c,16^a\}$, and~$L'' = \{1^c, 19^a,20^b\}$.  Any counterexample must have~$a<9$, $b<16$ and~$c<20$, with~$a+b+c=43-1=42$.  The only possibility is~$(a,b,c) = (8,15,19)$.  We obtain a realization for~$\{1^{19}, 19^8, 20^{15}\}$ from Lemma~\ref{lem:sawtooth}.

{\em Case~2:} $y=12$, $v=41$.  The three equivalent multisets are~$L = \{1^a,11^b,12^c\}$, $L' = \{ 1^b,15^a,16^c\}$, and~$L'' = \{1^c, 17^a,18^b\}$.  Any counterexample must have~$a<12$, $b<16$ and~$c<18$, with~$a+b+c=41-1=40$.  There are ten possibilities for~$(a,b,c)$:
$$(8,15,17),(9,14,17),(9,15,16),(10,13,17),(10,14,16),$$
\vspace{-9mm}
$$(10,15,15),(11,12,17),(11,13,16),(11,14,15),(11,15,14).$$  
In each case, $c \geq 17 - \min(a,b) + 1$ and we obtain a realization for~$L''$ from Lemma~\ref{lem:sawtooth}.

Note that it was sufficient to know the smallest and largest possible values for~$a$,~$b$ and~$c$; we did not need to know exactly which triples needed attention.   Where possible, we shall use this in the remaining cases.

{\em Case~3:} $y=12$, $v=53$.  The three equivalent multisets are~$L = \{1^a,11^b,12^c\}$, $L' = \{ 1^c,22^a,23^b\}$, and~$L'' = \{1^b, 23^c,24^a\}$.  As~$a+b+c=53-1=52$, any counterexample must have~$7 \leq a<12$, $19 \leq b<24$ and~$18 \leq c<23$. 
We have $b \geq 17 \geq  23 - \min(a,c) + 1$ and we obtain a realization for~$L''$ from Lemma~\ref{lem:sawtooth}.

{\em Case~4:} $y=12$, $v=79$.  The three equivalent multisets are~$L = \{1^a,11^b,12^c\}$, $L' = \{ 1^c,32^b,33^a\}$, and~$L'' = \{1^b, 36^a,37^c\}$.  As~$a+b+c=79-1=78$, any counterexample must have~$7 \leq a<12$, $19 \leq b<37$ and~$18 \leq c<33$. 
We have $b \geq 27 \geq  36 - \min(a,c) + 1$ and we obtain a realization for~$L''$ from Lemma~\ref{lem:sawtooth}.

{\em Case~5:} $y=13$, $v=41$.  The three equivalent multisets are~$L = \{1^a,12^b,13^c\}$, $L' = \{ 1^b,16^c,17^a\}$, and~$L'' = \{1^c, 18^b,19^a\}$.  As~$a+b+c=41-1=40$, any counterexample must have~$6 \leq a<13$, $10 \leq b<17$ and~$12 \leq c<19$. 
When~$c\geq 13$, we have $c \geq  18 - \min(a,b) + 1$ and we obtain a realization for~$L''$ from Lemma~\ref{lem:sawtooth}.  This leaves one subcase to consider: $(a,b,c) = (12, 16,12)$.  As~$12 \neq 16-1$, Lemma~\ref{lem:sawtooth} provides a realization for~$L''$ here too. 

{\em Case~6:} $y=13$, $v=67$.  The three equivalent multisets are~$L = \{1^a,12^b,13^c\}$, $L' = \{ 1^b,28^a,29^c\}$, and~$L'' = \{1^c, 30^b,31^a\}$.  As~$a+b+c=67-1=66$, any counterexample must have~$8 \leq a<13$, $25 \leq b<29$ and~$27 \leq c<31$. 
We have $c \geq 23 \geq  30 - \min(a,b) + 1$ and we obtain a realization for~$L''$ from Lemma~\ref{lem:sawtooth}.

{\em Case~7:} $y=13$, $v=89$.  The three equivalent multisets are~$L = \{1^a,12^b,13^c\}$, $L' = \{ 1^b,36^c,37^a\}$, and~$L'' = \{1^c, 41^a,42^b\}$.  As~$a+b+c=89-1=88$, any counterexample must have~$11 \leq a<13$, $35 \leq b<37$ and~$40 \leq c<42$. 
We have $c \geq 31 \geq  30 - \min(a,b) + 1$ and we obtain a realization for~$L''$ from Lemma~\ref{lem:sawtooth}.

{\em Case~8:} $y=14$, $v=67$. The three equivalent multisets are~$L = \{1^a,13^b,14^c\}$, $L' = \{ 1^c,23^b,24^a\}$, and~$L'' = \{1^b, 31^a,32^c\}$.  As~$a+b+c=67-1=66$, any counterexample must have~$11 \leq a<14$, $29 \leq b<32$ and~$22 \leq c< 24$. 
We have $b \geq 21 \geq  31 - \min(a,c) + 1$ and we obtain a realization for~$L''$ from Lemma~\ref{lem:sawtooth}.

{\em Case~9:} $y=15$, $v=73$. The three equivalent multisets are~$L = \{1^a,14^b,15^c\}$, $L' = \{ 1^b,25^c,26^a\}$, and~$L'' = \{1^c, 34^a,35^b \}$.  As~$a+b+c=73-1=72$, any counterexample must have~$13 \leq a<15$, $24 \leq b<26$ and~$33 \leq c< 35$. 
We have $c \geq 22 \geq  34 - \min(a,b) + 1$ and we obtain a realization for~$L''$ from Lemma~\ref{lem:sawtooth}.

{\em Case~10:} $y=16$, $v=41$.  This is equivalent to {\em Case~2}.

{\em Case~11:} $y=16$, $v=43$.  This is equivalent to {\em Case~1}.

{\em Case~12:} $y=16$, $v=73$. The three equivalent multisets are~$L = \{1^a,15^b,16^c\}$, $L' = \{ 1^c,31^b,32^a\}$, and~$L'' = \{1^b, 33^c,34^a \}$.  As~$a+b+c=73-1=72$, any counterexample must have~$8 \leq a<16$, $26 \leq b<34$ and~$24 \leq c< 32$. 
We have $b \geq 26 \geq  33 - \min(a,c) + 1$ and we obtain a realization for~$L''$ from Lemma~\ref{lem:sawtooth}.

{\em Case~13:} $y=16$, $v=103$.  The three equivalent multisets are~$L = \{1^a,15^b,16^c\}$, $L' = \{ 1^c,45^a,46^b\}$, and~$L'' = \{1^b, 47^c, 48^a \}$.  As~$a+b+c=103-1=102$, any counterexample must have~$10 \leq a<16$, $42 \leq b< 48$ and~$40 \leq c< 46$. 
We have $b \geq 38 \geq  47 - \min(a,c) + 1$ and we obtain a realization for~$L''$ from Lemma~\ref{lem:sawtooth}.

{\em Case~14:} $y=16$, $v=137$.  The three equivalent multisets are~$L = \{1^a,15^b,16^c\}$, $L' = \{ 1^c,59^b,60^a\}$, and~$L'' = \{1^b, 64^a, 65^c \}$.  As~$a+b+c=137-1=136$, any counterexample must have~$13 \leq a<16$, $62 \leq b< 65$ and~$57 \leq c< 60$. 
We have $b \geq 52 \geq  64 - \min(a,c) + 1$ and we obtain a realization for~$L''$ from Lemma~\ref{lem:sawtooth}.

{\em Case~15:} $y=16$, $v=167$. The three equivalent multisets are~$L = \{1^a,15^b,16^c\}$, $L' = \{ 1^c,73^a,74^b\}$, and~$L'' = \{1^b, 78^a, 79^c \}$.  As~$a+b+c=167-1=166$, as in Case~1, there is exactly one possibility for a counterexample, which is $(a,b,c) = (15,78,73)$ here.  Lemma~\ref{lem:sawtooth} gives a realization for~$L''$.
\end{proof}

This proves Theorem~\ref{th:small_y-1}: ``The Coprime BHR Conjecture holds for supports of the form~$\{1,y-1,y\}$ for all~$y \leq 16$."

\section{Extending the methods to the case of~$k=2$}\label{sec:y-2}

In this section, we prove that multisets of the form~$\{1^a, (y-2)^b, y^c\}$ have a linear realization for each~$a \geq y$.  The linear realizations we construct are standard when~$y$ is odd.  Coupled with the known constructions for even~$y$, this proves Theorem~\ref{th:biggie}.

The constructions follow a similar route to the~$k=1$ case.  We introduce certain moves, allowing us to switch between~$(y-2)$-edges and~$y$-edges, including a technique to switch between a higher number of them. (A similar technique was available for the~$k=1$ case as well, but was not needed.) We once again use these moves on~$\omega$-constructions, divided into separate lemmas based on the sizes of~$b$ and~$c$. Finally, we concatenate the modified~$\omega$-constructions with the appropriate standard linear realizations constructed in Section~\ref{sec:para} to get the desired results.

We switch between~($y-2$)-edges and~$y$-edges as follows: Suppose we have a subsequence~$(j-1, j, y+j-2, y+j-1, j+1)$ in a path, which realizes~$\{1^2, (y-2)^2\}$.  Without disturbing the rest of the path, we may replace this with the subsequence~$(j-1, y+j-1, y+j-2, j, j+1)$, which realizes~$\{1^2, y-2, y\}$.   Similarly, suppose we have a subsequence~$(j-1, j, y+j, y+j-1, j+1)$ in a path, which realizes~$\{1^2, y-2, y\}$.  Without disturbing the rest of the path, we may replace this with the subsequence~$(j-1, y+j-1, y+j, j, j+1)$, which realizes~$\{1^2, y^2\}$.  Each of these operations reverses a portion of the sequence.  We call such operations {\em $\beta$-moves}, consistent with the naming conventions of~\cite{AO1}.

The slanted grid-based representation (with a vertical step corresponding to a distance of~$y-1$) illustrates the geometry of these two $\beta$-moves that are horizontal reflections of one another:

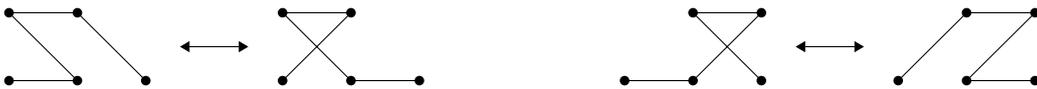
\begin{center}
\begin{tikzpicture}[scale=0.9, every node/.style={transform shape}]
\foreach \p in {1,5,10,14}{\fill (\p ,1) circle (2pt); \fill (\p+1 ,1) circle (2pt); \fill (\p+2 ,1) circle (2pt);  \fill (\p +1,2) circle (2pt);}
\foreach \p in {1,5,12,16}{ \fill (\p ,2) circle (2pt);}
\draw (1,1)--(2,1)--(1,2)--(2,2)--(3,1) ;
\draw (5,1)--(6,2)--(5,2)--(6,1)--(7,1) ;
\draw (10,1)--(11,1)--(12,2)--(11,2)--(12,1) ;
\draw (14,1)--(15,2)--(16,2)--(15,1)--(16,1);
\draw [Triangle-Triangle] (3.5,1.5) -- (4.5,1.5) ;
\draw [Triangle-Triangle] (12.5,1.5) -- (13.5,1.5) ;
\end{tikzpicture}
\end{center}

Lemma~\ref{lem:switch} introduces the technical widget that allows us to switch between higher number of~$(y-2)$-edges and~$y$-edges in linear realizations.

\begin{lem}\label{lem:switch} {\rm (Traversals of isosceles trapezoids)}
Let~$3\leq s<t$ with~$s$ odd.  For all~$b,c$ with~$b+c=s-1$ there is a path in~$K_v$ that uses the vertices $\{ 0,1, \ldots, s, t, t+1, \ldots, s+t-2\}$, that linearly realizes~$\{1^s, (t-2)^b, t^c\}$, and that has the following additional properties:
\begin{itemize}
\item the first vertex is $0$,
\item the final vertex is~$s$,
\item it includes the $1$-edges between $t+i$  and $t+i+1$ for even~$i$ in the range~$0 \leq i \leq s-3$. 
\end{itemize}
\end{lem}

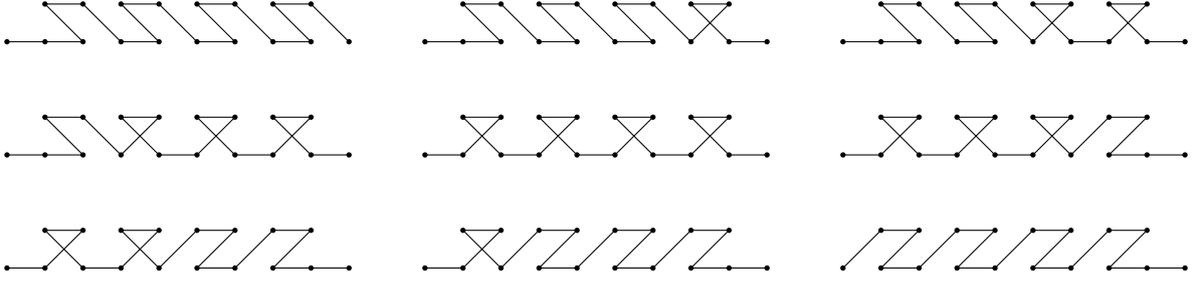
\begin{figure}
\caption{The path patterns realizing $\{1^9, (t-2)^b, t^c\}$ for $b+c=9$ from Lemma~\ref{lem:switch}.}\label{fig:switch9}
\begin{center}
\begin{tikzpicture}[scale=0.5, every node/.style={transform shape}]
\foreach \p in {1, 12, 23}{ \foreach \q in {1,4,7}{\foreach \r in {0,1,2,3,4,5,6,7,8,9}{\fill (\p + \r ,\q) circle (2pt); }}}
\foreach \p in {1, 12, 23}{ \foreach \q in {2,5,8}{\foreach \r in {1,2,3,4,5,6,7,8}{\fill (\p + \r ,\q) circle (2pt); }}}
\draw (1,7)--(3,7)--(2,8)--(3,8)--(4,7)--(5,7)--(4,8)--(5,8)--(6,7)--(7,7)--(6,8)--(7,8)--(8,7)--(9,7)--(8,8)--(9,8)--(10,7);
\draw (12,7)--(14,7)--(13,8)--(14,8)--(15,7)--(16,7)--(15,8)--(16,8)--(17,7)--(18,7)--(17,8)--(18,8)--(19,7)--(20,8)--(19,8)--(20,7)--(21,7);
\draw (23,7)--(25,7)--(24,8)--(25,8)--(26,7)--(27,7)--(26,8)--(27,8)--(28,7)--(29,8)--(28,8)--(29,7)--(30,7)--(31,8)--(30,8)--(31,7)--(32,7);
\draw (1,4)--(3,4)--(2,5)--(3,5)--(4,4)--(5,5)--(4,5)--(5,4)--(6,4)--(7,5)--(6,5)--(7,4)--(8,4)--(9,5)--(8,5)--(9,4)--(10,4);
\draw (12,4)--(13,4)--(14,5)--(13,5)--(14,4)--(15,4)--(16,5)--(15,5)--(16,4)--(17,4)--(18,5)--(17,5)--(18,4)--(19,4)--(20,5)--(19,5)--(20,4)--(21,4);
\draw (23,4)--(24,4)--(25,5)--(24,5)--(25,4)--(26,4)--(27,5)--(26,5)--(27,4)--(28,4)--(29,5)--(28,5)--(29,4)--(30,5)--(31,5)--(30,4)--(32,4);
\draw (1,1)--(2,1)--(3,2)--(2,2)--(3,1)--(4,1)--(5,2)--(4,2)--(5,1)--(6,2)--(7,2)--(6,1)--(7,1)--(8,2)--(9,2)--(8,1)--(10,1);
\draw (12,1)--(13,1)--(14,2)--(13,2)--(14,1)--(15,2)--(16,2)--(15,1)--(16,1)--(17,2)--(18,2)--(17,1)--(18,1)--(19,2)--(20,2)--(19,1)--(21,1);
\draw (23,1)--(24,2)--(25,2)--(24,1)--(25,1)--(26,2)--(27,2)--(26,1)--(27,1)--(28,2)--(29,2)--(28,1)--(29,1)--(30,2)--(31,2)--(30,1)--(32,1);
\end{tikzpicture}
\end{center}
\end{figure}

\begin{proof}
We start with the case~$c=0$ and successively tweak the construction, removing a $(y-2)$-edge and adding a $y$-edge each time.

For even~$i$ in the range~$0 \leq i \leq s-3$, define $\Theta_i$ to be the sequence $(i+1, i+2, t, t+i)$.  Then the path $(0, \Theta_0, \Theta_2, \ldots, \Theta_{s-3})$ satisfies the statement of the theorem with~$b=s-1$ and~$c=0$.  

We may apply the first $\beta$-move above up to~$(s-3)/2$ times.  Each instance replaces a~$\Theta_i$.   After~$(s-3)/2$ applications the path has the form $(\Theta'_0, \Theta'_2, \ldots, \Theta'_{s-3}, s)$, where $\Theta'_i = (i, i+1, y+i+1, y+i, i+2)$ for even~$i$ in the range $0 \leq i \leq s-3$.  We may now apply the second $\beta$-move a further~$(s-3)/2$ times, with each move replacing a~$\Theta'_i$.  

After~$j$ moves we have a path that linearly realizes the multiset $\{1^s, (t-2)^b, t^c\}$, where~$c = j$ and~$b = s-1-j$.  The moves do not affect the starting or ending vertices, or the $1$-edges between vertices labeled~$t$ or greater and so the additional claimed properties hold for all of these paths.

The nine paths for the case~$s=9$ are illustrated in Figure~\ref{fig:switch9}, where the lower row of vertices in each diagram is $\{0,1, \ldots 9\}$ and the upper row is $\{t, t+1, \ldots, t+7\}$.
\end{proof}

In the realizations used in the proof of Lemma~\ref{lem:switch}, the first edge is a 1-edge in all cases except one, and similarly, the final edge is a 1-edge in all cases but one.  The two extreme cases are those where~$b=0$ and where~$c=0$. Removing such a 1-edge gives Lemma~\ref{lem:switch2}, a useful variant of Lemma~\ref{lem:switch} for when there is a little less space to maneuver.

\begin{lem}\label{lem:switch2} {\rm (Traversals of right-trapezoids)}
Let~$3\leq s<t$ with~$s$ odd.  For all~$b,c > 0$ with~$b+c=s-1$ there is a path in~$K_v$ that uses the vertices $\{ 0,1, \ldots, s-1, t, t+1, \ldots, s+t-2\}$, that linearly realizes~$\{1^{s-1}, (t-2)^b, t^c\}$, and that has the following additional properties:
\begin{itemize}
\item the first vertex is $0$,
\item the final vertex is~$s-1$,
\item it includes the $1$-edges between $t+i$  and $t+i+1$ for even~$i$ in the range~$0 \leq i \leq s-3$. 
\end{itemize}
\end{lem}

We now introduce another diagrammatic tool to make the proofs of the main results of the section flow more smoothly. 
To indicate regions of a linear realization in which we can use the realizations given in Lemmas~\ref{lem:switch} and~\ref{lem:switch2} (or a translation thereof) we replace the explicit path with a dotted outline of the region, which will look like one of the following trapezoids (the first for Lemma~\ref{lem:switch} and the second and third for Lemma~\ref{lem:switch2}):

\begin{center}
\begin{tikzpicture}[scale=0.9, every node/.style={transform shape}]
\foreach \p in {0,...,15}{ \fill (\p ,0) circle (2pt);}
\foreach \p in {1, 7,11}{ \fill (\p ,1) circle (2pt); \fill (\p+1 ,1) circle (2pt);\fill (\p+2 ,1) circle (2pt);  \fill (\p+3 ,1) circle (2pt);}
\draw [thick, dotted] (0,0)--(5,0)--(4,1)--(1,1)--(0,0) ; 
\draw [thick, dotted] (6,0)--(10,0)--(10,1)--(7,1)--(6,0) ; 
\draw [thick, dotted] (11,0)--(15,0)--(14,1)--(11,1)--(11,0) ; 
\end{tikzpicture}
\end{center}

Using the notation of the lemmas, denote the path given by Lemma~\ref{lem:switch} by~$\Theta_{s,t,b,c}$ and the path given by Lemma~\ref{lem:switch2} by~$\Theta'_{s-1,t,b,c}$.  Such a shape in a grid-based graph representation of a linear realization means that we can fill that region in accordance with Lemmas~\ref{lem:switch} and~\ref{lem:switch2}.  In the diagrams in this section, we shall always use $t = y$ for the first two trapezoids and $t = y-1$ for the third.

Lemma~\ref{lem:small_bc} gives the first main ingredient of the proof of Theorem~\ref{th:biggie}.1 as an almost-immediate consequence of Lemma~\ref{lem:switch}.  It also illustrates the use of the trapezoid diagrams.

\begin{lem}\label{lem:small_bc} {\rm ($k=2$; small $b$ and $c$)}
Let $L = \{1^a, (y-2)^b, y^c\}$ with $b+c < y$.  Then~$L$ has a standard linear realization when~$a \geq y-1$.
\end{lem}

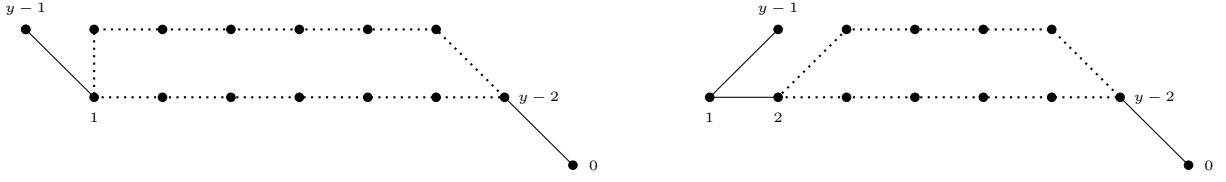
\begin{figure}
\caption{Diagrams for the subsidiary construction of Lemma~\ref{lem:small_bc}. }\label{fig:small_bc2}
\begin{center}
\begin{tikzpicture}[scale=0.9, every node/.style={transform shape}]
\foreach \p in {1,...,7}{ \fill (\p ,1) circle (2pt);}
\foreach \p in {0,...,6}{ \fill (\p ,2) circle (2pt);}
\foreach \p in {10,...,16}{ \fill (\p ,1) circle (2pt);}
\foreach \p in {11,...,15}{ \fill (\p ,2) circle (2pt);}
\fill (0,2) circle (2pt); \fill (8,0) circle (2pt); \fill (17,0) circle (2pt);
\draw [thick, dotted] (1,1)--(7,1)--(6,2)--(1,2)--(1,1) ; 
\draw [thick, dotted] (11,1)--(16,1)--(15,2)--(12,2)--(11,1) ; 
\draw (0,2)--(1,1); \draw (7,1)--(8,0);
\draw (11,2)--(10,1)--(11,1); \draw (16,1)--(17,0);
\node at (8.3, 0) {\tiny 0} ;   
\node at (0, 2.3) {\tiny $y-1$} ;   
\node at (1, 0.7) {\tiny $1$} ;   
\node at (7.5, 1) {\tiny $y-2$} ;  
\node at (10, 0.7) {\tiny 1} ;   
\node at (11, 2.3) {\tiny $y-1$} ;   
\node at (11, 0.7) {\tiny $2$} ;   
\node at (16.5, 1) {\tiny $y-2$} ;  
\node at (17.3, 0) {\tiny 0} ;   
\end{tikzpicture}
\end{center}
\end{figure}

\begin{proof}
It is sufficient to find a standard linear realization with~$a = y-1$ and apply Lemma~\ref{lem:concat1}.  If~$b$ or~$c$ is~0 then the $\omega$-construction of Lemma~\ref{lem:omega1} is sufficient, so assume~$b,c > 0$.

The main construction varies slightly with the parity of~$b+c$, with a subsidiary construction needed for the largest allowable value of $b+c$ when~$y$ is odd.  

{\em Case 1: $b + c$ even, $b+c \neq y-1$}.  The Hamiltonian path 
$$(\Theta_{b+c+1 , y, b,c}, b+c+1, b+c+2, \ldots, y-1),$$
which is the construction of Lemma~\ref{lem:switch} giving the correct number of $(y-2)$ and $y$-edges, followed by a series of 1-edges is a standard linear realization for~$L$.

{\em Case 2: $b + c$ odd}.  We take the path~$\Theta_{b+c,y, b, c-1}$, which gives all of the required $(y-2)$-edges and all but one of the required~$y$-edges, and translate it up by~1.  Append 1-edges as necessary as in the even case and attach the vertex~$y$ to~0; that is:
$$(\Theta_{b+c , y, b,c-1} + 1, b+c, b+c+1, b+c+2, \ldots, y, 0).$$
Reversing this gives the required standard linear realization for~$L$.

Diagrams illustrating the two variations of the main construction are given in Figure~\ref{fig:small_bc}.

\begin{figure}
\caption{Diagrams for the main construction of Lemma~\ref{lem:small_bc}. }\label{fig:small_bc}
\begin{center}
\begin{tikzpicture}[scale=0.9, every node/.style={transform shape}]
\foreach \p in {0,...,7}{ \fill (\p ,1) circle (2pt);}
\foreach \p in {1,...,4}{ \fill (\p ,2) circle (2pt);}
\foreach \p in {9,...,16}{ \fill (\p ,1) circle (2pt);}
\foreach \p in {10,...,13}{ \fill (\p ,2) circle (2pt);}
\fill (15,0) circle (2pt);
\draw [thick, dotted] (0,1)--(5,1)--(4,2)--(1,2)--(0,1) ; 
\draw [thick, dotted] (9,1)--(14,1)--(13,2)--(10,2)--(9,1) ; 
\draw (5,1)--(7,1);
\draw (14,1)--(16,1)--(15,0);
\node at (0, 0.7) {\tiny 0} ;   
\node at (1, 2.3) {\tiny $y$} ;   
\node at (5, 0.7) {\tiny $b+c$} ;   
\node at (7, 0.7) {\tiny $y-1$} ;  
\node at (9, 0.7) {\tiny 1} ;   
\node at (10, 2.3) {\tiny $y+1$} ;   
\node at (14, 0.7) {\tiny $b+c$} ;   
\node at (16, 1.3) {\tiny $y$} ;  
\node at (14.7, 0) {\tiny 0} ;   
\end{tikzpicture}
\end{center}
\end{figure}
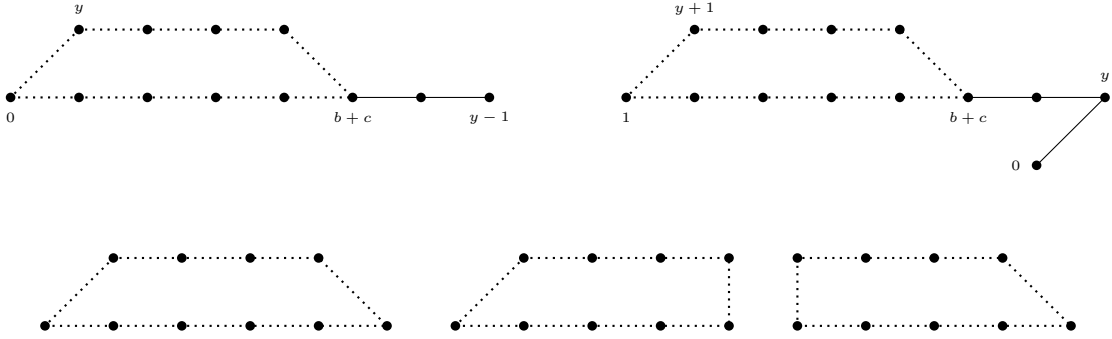

{\em Case 3: $b+c = y-1$ for odd~$y$}.  The construction depends on $c$.

First suppose~$c = 1$, and so~$b = y-2$.  Take the path~$\Theta'_{b-1,y, b-2, 1}$, which realizes the required $y$-edge and all but two of the required $(y-2)$-edges, and translate it up by~1.  Add~$(y-2)$-edges from~0 to~$y-2$ and from~1 to~$y-1$ to give the path
$$(y-1, \Theta'_{b-1,y, b-2, 1} + 1 , 0)$$
and reverse to give the required standard linear realization.

Now suppose~$c>1$.   Take the path~$\Theta_{b+c-2,y-1, b, c-2}$, which realizes  the required $(y-2)$-edges and all but two of the required $y$-edges, and translate it up by~2.  Add~$y$-edges from~0 to~$y$ and from~1 to~$y+1$, and a 1-edge from~1 to~2 to give the path
$$(y+1, 1, 2, \Theta_{b-1,y-1, b-2, 1} + 2 , 0)$$
and reverse to give the required standard linear realization.

Diagrams for these constructions are given in Figure~\ref{fig:small_bc2}.
\end{proof}

We can now move on to the main constructions.

\begin{lem}\label{lem:small_b} {\rm ($k=2$; small $b$)}
Let $L = \{1^a, (y-2)^b, y^c\}$ with,~$b \leq y-1$ and~$b+c \geq y$.  Then~$L$ has a standard linear realization when $a \geq y$.
\end{lem}

\begin{proof}
If either~$b$ or~$c$ is~0 we may use the $\omega$-constructions of Section~\ref{sec:background}, so assume that $b, c > 0$.  Using Lemma~\ref{lem:concat1}, it is sufficient to find standard linear realizations for any value of~$a$ with~$a \leq y$; in some cases, the constructions have~$a$ strictly less than~$y$ giving a slightly stronger result when the relevant conditions apply.
Similarly, sometimes the constructions naturally cover cases with~$b > y-1$.

Write~$b+c = qy + r$ with $0 \leq r < y$.  The conditions imply that~$q \geq 1$.   The proof splits into four cases depending on the parities of~$y$ and~$r$.

{\em Case 1: $y$ and $r$ even.}   Suppose~$y>4$. We construct standard linear realizations of type~$\mathcal{C}_{y}$ when~$b+c = y$ for all values of~$b$ in the range~$1 \leq b \leq y-1$. This case requires two constructions.   The available parameter ranges for the two constructions overlap, so there is choice of which to use for~$b$ near the middle of the ranges.

First, consider~$1 \leq b \leq y-3$ (and so~$c = y-b$ is in the range $3 \leq c \leq y-1$).  Let~$v = 2y$.  Take the path~$\Theta'_{y-2,y,b,c-2}$ from Lemma~\ref{lem:switch}, which has last vertex $y-2$, realizes~$\{1^{y-2}, (y-2)^b, y^{c-2}\}$ and omits the vertices~$\{y-1, 2y-2, 2y-1\}$.  
The path~$(\Theta'_{y-2,y,b,c-2}, 2y-2, 2y-1, y-1)$ is Hamiltonian and realizes $$\{1^{y-2}, (y-2)^b, y^{c-2}\} \cup \{1, y^2\} = \{1^{y-1}, (y-2)^b, y^{c}\}$$
as required.  Further, it is of type~$\mathcal{C}_y$, by the properties guaranteed in Lemma~\ref{lem:switch2} and the inclusion of the edge~$(2y-2, 2y-1)$.  

Second, consider~$4 \leq b \leq y-1$ (and so~$c = y-b$ is in the range $1 \leq c \leq y-4$).  Let~$v = 2y$.
Take the path~$\Theta'_{y-4,y-1,b-3,c-1}$ from Lemma~\ref{lem:switch} and consider $\Theta'_{y-4,y-1,b-3,c-1} + 1$,
which has last vertex $y-3$, realizes~$\{1^{y-4}, (y-2)^{b-3}, y^{c-1}\}$ and omits the vertices
$$\{0,y-2, y-1, 2y-4, 2y-3, 2y-2, 2y-1\}.$$  
The path
$$(0, y-2, 2y-4, 2y-3, 2y-2, 2y-1, y-1,  \Theta'_{y-4,y-1,b-3,c-1} + 1)$$ 
is Hamiltonian and realizes 
$$\{1^3, (y-2)^3, y \} \cup  \{1^{y-4}, (y-2)^{b-3}, y^{c-1}\}  = \{1^{y-1}, (y-2)^b, y^{c}\}$$
as required. Further, it is of type~$\mathcal{C}_y$, by the properties guaranteed in Lemma~\ref{lem:switch2} and the inclusion of the edges $(2y-4, 2y-3)$ and~$(2y-2, 2y-1)$.  

The two constructions are illustrated in  Figure~\ref{fig:smb_ev_ev}.  The diagrams use~$y=10$.  The diagrams for arbitrary~$y$ are the same, but with suitably elongated (or shortened) trapezoids.

\begin{figure}
\caption{Standard linear realizations of type $\mathcal{C}_y$ for {\em Case~1} of Lemma~\ref{lem:small_b}. }\label{fig:smb_ev_ev}
\begin{center}
\begin{tikzpicture}[scale=0.7, every node/.style={transform shape}]
\foreach \p in {0,...,9}{ \fill (\p ,0) circle (2pt);}
\foreach \p in {14,...,21}{ \fill (\p ,0) circle (2pt);}
\foreach \p in {1,...,10}{ \fill (\p ,1) circle (2pt);}
\foreach \p in {13,...,20}{ \fill (\p ,1) circle (2pt);}
\foreach \p in {12,...,15}{ \fill (\p ,2) circle (2pt);}
\draw [thick, dotted] (0,0)--(8,0)--(8,1)--(1,1)--(0,0) ; 
\draw [thick, dotted] (15,0)--(21,0)--(20,1)--(15,1)--(15,0) ; 
\draw (8,0) -- (9,1) -- (10,1)--(9,0);
\foreach \p in {1,3,...,9} {\draw [very thick] (\p,1)--(\p+1,1) ; \draw [-stealth]  (\p+0.5,1) -- (\p+0.8, 1.3) ;}
\draw (14,0) -- (12,2) -- (15,2) -- (14,1) -- (15,0) ;
\foreach \p in {12,14} {\draw [very thick] (\p,2)--(\p+1,2) ; \draw [-stealth]  (\p+0.5,2) -- (\p+0.8, 2.3) ;}
\foreach \p in {15,17,19} {\draw [very thick] (\p,1)--(\p+1,1) ; \draw [-stealth]  (\p+0.5,1) -- (\p+0.8, 1.3) ;}
\end{tikzpicture}
\end{center}
\end{figure}
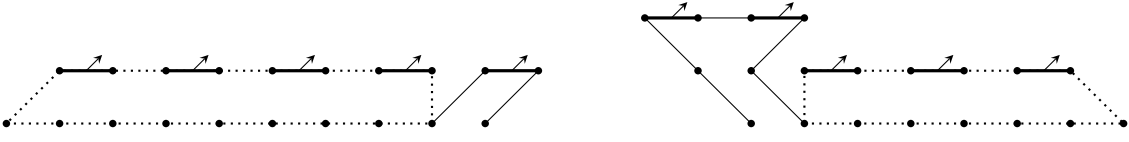

The first diagram of Figure~\ref{fig:smb_small} covers the subcase~$y=4$ in the same way.

\begin{figure}
\caption{Standard linear realizations with small~$y$ for {\em Cases~1-2} and~{\em 4a} of Lemma~\ref{lem:small_b}. }\label{fig:smb_small}
\begin{center}
\begin{tikzpicture}[scale=0.7, every node/.style={transform shape}]
\foreach \p in {8,13,18,19,20}{ \fill (\p ,0) circle (2pt);}
\foreach \p in {1,2,6,7,10,11,12,13,14,17,18,19,20,21}{ \fill (\p ,1) circle (2pt);}
\foreach \p in {0,1,2,3,5,6,7,8,11,12,13,14,15,18,19,20}{ \fill (\p ,2) circle (2pt);}
\foreach \p in {1,2,6,7}{ \fill (\p ,3) circle (2pt);}
\draw (1,1)--(0,2)--(1,2)--(2,3)--(1,3)--(2,2)--(3,2)--(2,1);
\draw  (6,1)--(5,2)--(6,2)--(7,3)--(6,3)--(7,2)--(8,2)--(7,1)--(8,0);
\draw (11,2)--(10,1)--(11,1)--(12,2)--(13,2)--(14,1)--(15,2)--(14,2)--(13,1)--(12,1)--(13,0);
\draw (18,0)--(17,1)--(18,1)--(19,0)--(20,1)--(21,1)--(20,0)--(18,2)--(20,2);
 \foreach \p in {1,6} {\draw [very thick] (\p,3)--(\p+1,3) ; \draw [-stealth]  (\p+0.5,3) -- (\p+0.8, 3.3) ;  \draw [-stealth]  (\p+1,3) -- (\p+1.3, 3.3) ; \draw (\p,3) circle (4pt);  \draw (\p+1,3) circle (4pt); }
 \foreach \p in {2,7} {\draw [very thick] (\p,2)--(\p+1,2) ; \draw [-stealth]  (\p+0.5,2) -- (\p+0.8, 2.3) ;  \draw (\p,2) circle (4pt);  \draw (\p+1,2) circle (4pt);} 
\foreach \p in {11,...,15} { \draw (\p,2) circle (4pt);} \foreach \p in {18,...,20} { \draw (\p,2) circle (4pt);} \foreach \p in {20,21} { \draw (\p,1) circle (4pt);}
\draw  [-stealth] (15,2)--(15.3,2.3);  \draw  [-stealth] (20,2)--(20.3,2.3); 
\end{tikzpicture}
\end{center}
\end{figure}
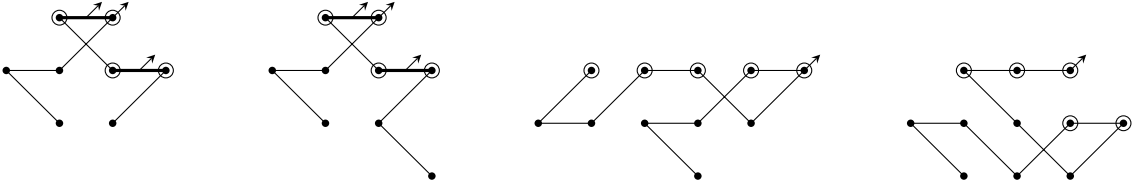

{\em Case 2: $y$ even and $r$ odd.}  Note that the parity conditions imply we in fact have $b+c \geq y+1$.  We adapt the constructions from {\em Case~1} into standard linear realizations of type $\mathcal{C}_y$ for~$b+c = y+1$.  Let~$v = 2y+1$.

The adaptation is the same in each subcase, with slightly different outcomes: translate the linear realization up by~1 and add an edge from the final vertex to~0.

In the first subcase, with~$y>4$ and~$1 \leq b \leq y-3$, we obtain
$(\Theta'_{y-2,y,b,c-2}+1, 2y-1, 2y, y, 0),$
which adds a $y$-edge and so realizes $\{1^{y-1}, (y-2)^b, y^{c}\}$
for~$b$ in the same range (and with~$c$ one larger than previously).

In the second subcase, the process adds a $(y-2)$-edge.  The path is 
$$(1, y-1, 2y-3, 2y-2, 2y-1, 2y, y,  \Theta'_{y-4,y-1,b-3,c-1} + 2, 0)$$ 
which realizes~$\{1^{y-1}, (y-2)^b, y^{c}\}$ for~$5 \leq b \leq y$.
Reverse these paths to make them standard (without changing the realized multiset or their type~$\mathcal{C}_y$ status).

Finally, the second diagram of Figure~\ref{fig:smb_small} covers the subcase~$y=4$.

{\em Case 3: $y$ odd and $r$ even.}  Similarly to the previous cases, we have two constructions that between them create standard linear realizations of type~$\mathcal{C}_y$ with~$b+c = y$ for~$b$ in the range~$1 \leq b \leq y-1$.  There is again overlap near the middle of the range.

First, consider~$1 \leq b \leq y-2$ and let~$v = 2y$.  Take the path $\Theta'_{y-1,y,b,c-1}$, which has last vertex~$y$, realizes~$\{1^{y-1}, (y-2)^b, y^{c-1}\}$ and omits the vertex~$2y-1$.  The path~$( \Theta'_{y-1,y,b,c-1}, 2y-1)$ is Hamiltonian and realizes   
$$\{1^{y-1}, (y-2)^b, y^{c-1}\} \cup \{ y \} = \{1^{y-1}, (y-2)^b, y^{c}\}$$
as required.  It is of type~$\mathcal{C}_y$ by the properties of~$\Theta'_{y-1,y,b,c-1}$ guaranteed by Lemma~\ref{lem:switch2} and that the final vertex is~$2y-1 = v-1$.

Second, consider~$4 \leq b \leq y-1$ and let~$v=2y$. Take the path~$\Theta'_{y-3,y,b-3,c}$ from Lemma~\ref{lem:switch2} and consider $\Theta'_{y-3,y,b-3,c} + 1$, which has last vertex~$y-2$, realizes $\{1^{y-3}, (y-2)^{b-3}, y^c\}$ and omits the vertices
$$\{0, y-1, 2y-3,2y-2,2y-1)$$
The path
$$( 0, \Theta'_{y-3,y,b-3,c} + 1, y-1, 2y-3,2y-2,2y-1)$$
is Hamiltonian and realizes
$$\{y-2\} \cup \{1^{y-3}, (y-2)^{b-3}, y^c\} \cup \{1^2, (y-2)^2\} = \{1^{y-1}, (y-2)^{b}, y^c\}$$
as required.  Further, it is of type~$\mathcal{C}_y$ by Lemma~\ref{lem:switch2}, the inclusion of the edge~$(2y-3,2y-2)$, and having final vertex~$2y-1 = v-1$.

The two constructions for this case are illustrated in Figure~\ref{fig:smb_odd_ev}.  The diagrams use~$y=9$; suitably elongated (or shortened) trapezoids transform them to the constructions for arbitrary $y$.

\begin{figure}
\caption{Standard linear realizations of type $\mathcal{C}_y$ for {\em Case~3} of Lemma~\ref{lem:small_b}. }\label{fig:smb_odd_ev}
\begin{center}
\begin{tikzpicture}[scale=0.7, every node/.style={transform shape}]
\fill (20,0) circle (2pt);
\foreach \p in {0,...,8}{ \fill (\p ,1) circle (2pt);}
\foreach \p in {13,...,19}{ \fill (\p ,1) circle (2pt);}
\foreach \p in {1,...,9}{ \fill (\p ,2) circle (2pt);}
\foreach \p in {12,...,18}{ \fill (\p ,2) circle (2pt);}
\foreach \p in {11,...,13}{ \fill (\p ,3) circle (2pt);}
\draw [thick, dotted] (0,1)--(8,1)--(8,2)--(1,2)--(0,1) ; 
\draw [thick, dotted] (13,1)--(19,1)--(18,2)--(13,2)--(13,1) ; 
\draw (8,1) -- (9,2);
\draw (20,0) -- (19,1);
\draw (13,1) -- (12,2) -- (11,3) -- (13,3) ;
\foreach \p in {1,3,...,7} {\draw [very thick] (\p,2)--(\p+1,2) ; \draw [-stealth]  (\p+0.5,2) -- (\p+0.8, 2.3) ;}
 \draw [-stealth]  (9,2) -- (9.3, 2.3) ;
\draw [very thick] (11,3)--(12,3) ; \draw [-stealth]  (11.5,3) -- (11.8, 3.3) ;
 \draw [-stealth]  (13,3) -- (13.3, 3.3) ;
 \foreach \p in {13,15,17} {\draw [very thick] (\p,2)--(\p+1,2) ; \draw [-stealth]  (\p+0.5,2) -- (\p+0.8, 2.3) ;}
\end{tikzpicture}
\end{center}
\end{figure}
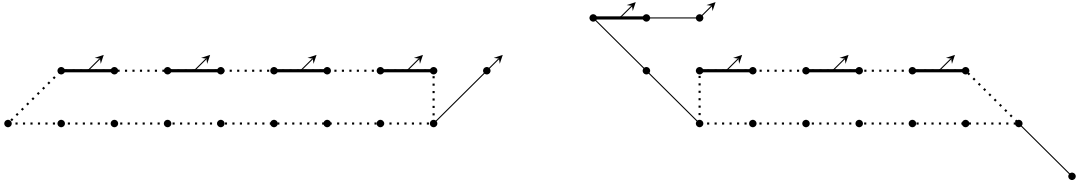

{\em Case 4a: $y$ odd and~$r=1$.}  We construct $y$-growable standard linear realizations when~$b+c = y+1$.  There are four main subcases, the first two of which overlap.

First, take $y>5$ and $1 \leq b \leq y-4$.  We use a translation of $\Theta'_{y-3,y,b, c-4}$ with some surrounding edges:
$$(y+1, 1, \Theta'_{y-3,y,b, c-4} + 2, 2y-1, 2y, y, 0).$$
It realizes 
$$\{1,y\} \cup \{1^{y-3}, (y-2)^b, y^{c-4}\} \cup
 \{ 1, y^3 \} = \{1^{y-1}, (y-2)^b, y^{c}\}$$
and reversing gives a standard linear realization.

Second, take $y>5$ and $3 \leq b \leq y-2$.  For this subcase we need to remove a 1-edge from each end of the sequence from Lemma~\ref{lem:switch}, rather than just one end or other as in Lemma~\ref{lem:switch2}.  Diagrammatically, this makes our trapezoid a rectangle; extend the $\Theta$ notation to this situation in the obvious way, using $\Theta''_{s-2, t, b, c}$ to indicate the path with first vertex~0, last vertex~$s-2$ and realizing $\{1^{s-2}, (t-2)^b, t^c\}$ when~$b,c > 0$ and $b+c = s-1$ for odd~$s$. 
Then, consider the Hamiltonian path
$$(y-1, \Theta''_{y-4,y-1,b-2,c-2}, 2y-3, 2y-2, y-2, 0)$$
in~$K_{2y-1}$.
This realizes
$$\{y-2\} \cup \{1^{y-4}, (y-2)^{b-2}, y^{c-2}\} \cup \{1, y-2, y^2\} = \{ 1^{y-3}, (y-2)^b, y^c\}$$
and the reverse is standard.

The next two subcases use a specific instance of Lemma~\ref{lem:switch2} (using a single $y$-edge in each subcase) for which it is straightforward to confirm the required properties.
When~$y>5$ and $b = y-1$, take the reverse of 
$$(y,2,1,y-1,2y-3,2y-2,2y-1,2y,2y+1,y+1,3, 
\Theta'_{y-5,y-2, y-4,1}, 0),$$
which realizes~$\{1^y, (y-2)^b, y^c\}$.
When $y>5$ and $b = y$, take the reverse of
$$(y,2,1,y-1,2y-3,2y-2,2y-1,y+1,3, 
\Theta'_{y-5,y-2, y-4,1}, 0),$$
which realizes~$\{1^{y-2}, (y-2)^b, y^c\}$.

The four main constructions are illustrated in Figure~\ref{fig:smb_4a}.   The necessary constructions for~$y=5$ are given in the last two diagrams of Figure~\ref{fig:smb_small}.

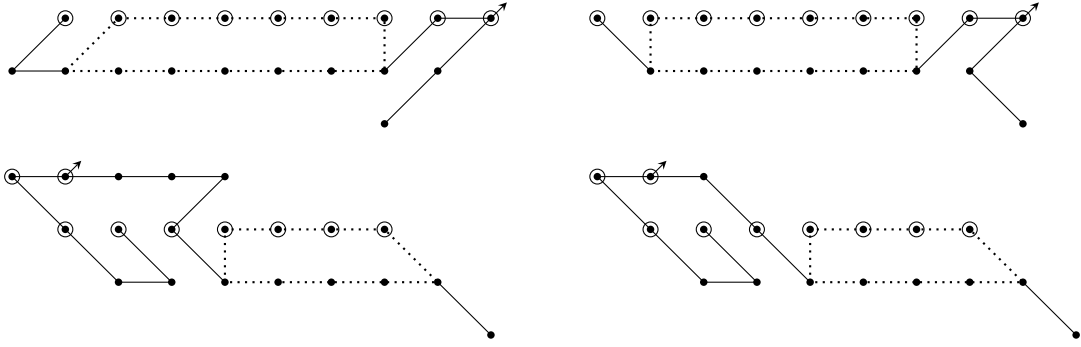
\begin{figure}
\caption{$y$-Growable standard linear realizations for {\em Case~4a} of Lemma~\ref{lem:small_b}. }\label{fig:smb_4a}
\begin{center}
\begin{tikzpicture}[scale=0.7, every node/.style={transform shape}]
\foreach \p in {0,...,8}{ \fill (\p ,5) circle (2pt);}
\foreach \p in {1,...,9}{ \fill (\p ,6) circle (2pt);}
\foreach \p in {12,...,18}{ \fill (\p ,5) circle (2pt);}
\foreach \p in {11,...,19}{ \fill (\p ,6) circle (2pt);}
\fill (7,4) circle (2pt); \fill (19,4) circle (2pt);
\draw [thick, dotted] (1,5)--(7,5)--(7,6)--(2,6)--(1,5) ; 
\draw [thick, dotted] (12,5)--(17,5)--(17,6)--(12,6)--(12,5) ; 
\draw (1,6)--(0,5)--(1,5); \draw (7,5)--(8,6)--(9,6)--(7,4);
\draw (11,6)--(12,5); \draw (17,5)--(18,6)--(19,6)--(18,5)--(19,4);
\foreach \p in {2,...,8}{ \fill (\p ,1) circle (2pt);}
\foreach \p in {1,...,7}{ \fill (\p ,2) circle (2pt);}
\foreach \p in {0,...,4}{ \fill (\p ,3) circle (2pt);}
\foreach \p in {13,...,19}{ \fill (\p ,1) circle (2pt);}
\foreach \p in {12,...,18}{ \fill (\p ,2) circle (2pt);}
\foreach \p in {11,...,13}{ \fill (\p ,3) circle (2pt);}
\fill (9,0) circle (2pt); \fill (20,0) circle (2pt);
\draw (2,2)--(3,1)--(2,1)--(0,3)--(4,3)--(3,2)--(4,1); \draw (8,1)--(9,0);
\draw (13,2)--(14,1)--(13,1)--(11,3)--(13,3)--(15,1); \draw (19,1)--(20,0);
\draw [thick, dotted] (4,1)--(4,2)--(7,2)--(8,1)--(4,1) ; 
\draw [thick, dotted] (15,1)--(15,2)--(18,2)--(19,1)--(15,1) ; 

\foreach \p in {1,...,9}{\draw (\p,6) circle (4pt);}
\draw [-stealth] (9,6)--(9.3,6.3);
\foreach \p in {11,...,19}{\draw (\p,6) circle (4pt);}
\draw [-stealth] (19,6)--(19.3,6.3);

\foreach \p in {0,1}{\draw (\p,3) circle (4pt);}
\foreach \p in {1,...,7}{\draw (\p,2) circle (4pt);}
\draw [-stealth] (1,3)--(1.3,3.3);

\foreach \p in {11,12}{\draw (\p,3) circle (4pt);}
\foreach \p in {12,...,18}{\draw (\p,2) circle (4pt);}
\draw [-stealth] (12,3)--(12.3,3.3);

\end{tikzpicture}
\end{center}
\end{figure}

{\em Case 4b: $y$ odd and $r >1$.}  
This subcase has four main constructions that cover $y>7$ and then variations to cover~$y \in \{5,7\}$. 
For the main constructions, we first construct $y$-growable standard linear realizations for~$b+c = y+3$, the smallest required value to cover, which has $r=3$.  We then show how to transform these constructions to reach arbitrary odd~$r$ in the range $3 < r < y$.

First, suppose $y>7$ and $4 \leq b \leq y-6$.  Start with the path~$\Theta'_{y-5,y,b,c-8}$, realizing~$\{1^{y-5}, (y-2)^b, y^{c-8} \}$.  We translate this up by~4 and use it in the following path:
$$(y+2,2,1,y+1,2y+1,2y+2,2y+3,y+3,3, \Theta'_{y-5,y,b,c-8}+4, 2y-1, 2y,y,0).$$
The reverse is standard, $y$-growable at~$2y$, and realizes
$$\{1^4, y^5\} \cup \{1^{y-5}, (y-2)^{b}, y^{c-8}\} \cup \{1, y^3\} = \{ 1^{y}, (y-2)^b, y^c\}.$$

Second, suppose $y>7$ and $4 \leq b \leq y-4$.  Start with the path~$\Theta'_{y-7,y-1,b-3,c-7}$, realizing~$\{1^{y-7}, (y-2)^{b-3}, y^{c-7} \}$.  We translate this up by~5 and use it in the following path:
$$(y+2,4,3,2,1,y+1,2y+1,2y+2,2y+3,y+3,\Theta'_{y-7,y-1,b-3,c-7}+5,$$
$$2y-2,2y-3,y-1,2y-1,2y,y,0).$$
The reverse is standard, $y$-growable at~$2y$, and realizes
$$\{1^5,(y-2)^2, y^3\} \cup \{1^{y-7}, (y-2)^{b-3}, y^{c-7}\} \cup \{1^2,y-2, y^4\} = \{ 1^{y}, (y-2)^b, y^c\}.$$

Third, suppose $y>7$ and $6 \leq b \leq y-2$.  Start with the path~$\Theta'_{y-7,y-1,b-5,c-5}$, realizing~$\{1^{y-7}, (y-2)^{b-5}, y^{c-5} \}$.  We translate this up by~3 and use it in the following path:
$$(y,2,1,y-1,2y-1,2y,2y+1,y+1,\Theta'_{y-7,y-1,b-4,c-5}+3,$$
$$2y-4,2y-5,y-3,2y-3,2y-2,y-2,0).$$
The reverse is standard, $y$-growable at~$2y-2$, and realizes
$$\{1^3,(y-2)^3, y^2\} \cup \{1^{y-7}, (y-2)^{b-5}, y^{c-5}\} \cup \{1^2,(y-2)^2, y^3\} = \{ 1^{y-2}, (y-2)^b, y^c\}.$$

Fourth, suppose $y>7$ and $8 \leq b \leq y$.  Start with the path~$\Theta'_{y-7,y-1,b-7,c-3}$, realizing~$\{1^{y-7}, (y-2)^{b-7}, y^{c-3} \}$.  We translate this up by~5 and use it in the following path:
$$(y,2,1,y-1,2y-3,2y-2,2y-1,2y,2y+1,y+1,3,4,y+2,2y+2,2y+3,y+3,$$
$$\Theta'_{y-7,y-1,b-7,c-3}, 0
).$$
The reverse is standard, $y$-growable at~$2y-1$, and realizes
$$\{1^7,(y-2)^6, y^3\} \cup \{1^{y-7}, (y-2)^{b-7}, y^{c-3}\} \cup \{y-2\} = \{ 1^{y}, (y-2)^b, y^c\}.$$

The four constructions are illustrated in Figure~\ref{fig:smb_4b}.

\begin{figure}
\caption{$y$-Growable standard linear realizations for {\em Case~4b} of Lemma~\ref{lem:small_b}.
}\label{fig:smb_4b}
\begin{center}
\begin{tikzpicture}[scale=0.6, every node/.style={transform shape}]
\foreach \p in {2,...,10}{ \fill (\p ,1) circle (2pt);}
\foreach \p in {15,...,23}{ \fill (\p ,1) circle (2pt);}
\foreach \p in {1,...,11}{ \fill (\p ,2) circle (2pt);}
\foreach \p in {14,...,22}{ \fill (\p ,2) circle (2pt);}
\foreach \p in {2,...,4}{ \fill (\p ,3) circle (2pt);}
\foreach \p in {13,...,19}{ \fill (\p ,3) circle (2pt);}
\foreach \p in {0,...,10}{ \fill (\p ,5) circle (2pt);}
\foreach \p in {13,...,23}{ \fill (\p ,5) circle (2pt);}
\foreach \p in {1,...,11}{ \fill (\p ,6) circle (2pt);}
\foreach \p in {14,...,24}{ \fill (\p ,6) circle (2pt);}
\foreach \p in {2,...,4}{ \fill (\p ,7) circle (2pt);}
\foreach \p in {15,...,17}{ \fill (\p ,7) circle (2pt);}
\fill (9,4) circle (2pt);\fill (22,4) circle (2pt);\fill (11,0) circle (2pt);\fill (24,0) circle (2pt);
\draw (2,6)--(1,5)--(0,5)--(2,7)--(4,7)--(2,5)--(3,5); \draw (9,5)--(10,6)--(11,6)--(9,4);
\draw (15,6)--(16,5)--(13,5)--(15,7)--(17,7)--(16,6)--(17,5); \draw (21,5)--(22,6)--(21,6)--(22,5)--(23,6)--(24,6)--(22,4);
\draw (2,2)--(3,1)--(2,1)--(1,2)--(2,3)--(4,3)--(3,2)--(4,1); \draw (8,1)--(9,2)--(8,2)--(9,1)--(10,2)--(11,2)--(10,1)--(11,0);
\draw (15,2)--(16,1)--(15,1)--(13,3)--(17,3)--(16,2)--(17,1)--(18,1)--(17,2)--(18,3)--(19,3)--(18,2)--(19,1); \draw (23,1)--(24,0);
\draw [thick, dotted]  (3,5)--(4,6)--(9,6)--(9,5)--(3,5);
\draw [thick, dotted]  (17,5)--(17,6)--(20,6)--(21,5)--(17,5);
\draw [thick, dotted]  (4,1)--(4,2)--(7,2)--(8,1)--(4,1);
\draw [thick, dotted] (19,1)--(19,2)--(22,2)--(23,1)--(19,1);
\foreach \p in {1,...,11}{\draw (\p,6) circle (4pt); \draw (\p,2) circle (4pt);}
\foreach \p in {14,...,24}{\draw (\p,6) circle (4pt);}
\foreach \p in {13,...,15}{\draw (\p,3) circle (4pt);} \foreach \p in {15,...,22}{\draw (\p,2) circle (4pt);}
 \foreach \p in {4,6,8,17,19,21} {\draw [very thick] (\p,6)--(\p+1,6) ; \draw [-stealth]  (\p+0.5,6) -- (\p+0.8, 6.3) ;}
 \foreach \p in {4,6,8} {\draw [very thick] (\p,2)--(\p+1,2) ; \draw [-stealth]  (\p+0.5,2) -- (\p+0.8, 2.3) ;}
\draw [very thick] (13,3)--(14,3) ; \draw [-stealth]  (13.5,3) -- (13.8, 3.3) ; \draw [very thick] (19,2)--(20,2) ; \draw [-stealth]  (19.5,2) -- (19.8, 2.3) ;  \draw [very thick] (21,2)--(22,2) ; \draw [-stealth]  (21.5,2) -- (21.8, 2.3) ;
 \draw [-stealth]  (11,6) -- (11.3, 6.3) ;   \draw [-stealth]  (24,6) -- (24.3, 6.3) ;  \draw [-stealth]  (11,2) -- (11.3, 2.3) ;  \draw [-stealth]  (15,3) -- (15.3, 3.3) ;

\end{tikzpicture}
\end{center}
\end{figure}
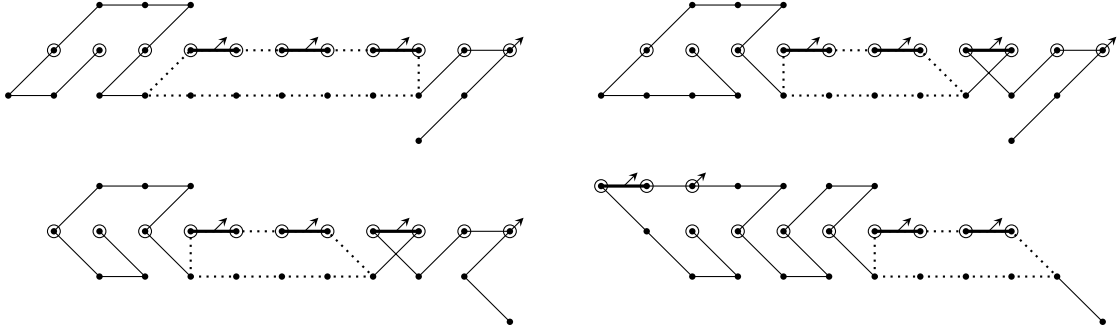

In the diagrams, notice the bold lines and accompanying arrows, in the style of type~$\mathcal{C}_y$ linear realizations, although these are not, in general, of type~$\mathcal{C}_y$.  However, they function in a similar way locally at the indicated edges.  

Specifically, each construction includes each of the edges $(v-y+i, v-y+i+1)$ for even~$i$ in the range $0 \leq i \leq y-7$.  We may replace each of these edges, starting with~$i=0$ and working through increasing values in turn, with the three-edge sequence $(v-y+i, v+i, v+i+1, v-y+i+1)$.   Each such replacement adds two $y$-edges to the linear realization, covering the required values of~$r$ with~$r>3$ in turn. 

\begin{figure}
\caption{Standard linear realizations with small~$y$ for {\em Case~4b} of Lemma~\ref{lem:small_b}. }\label{fig:smb_small2}
\begin{center}
\begin{tikzpicture}[scale=0.67, every node/.style={transform shape}]
\foreach \p in {0,...,4}{ \fill (\p ,0) circle (2pt);  \draw (\p,0) circle (4pt);}
\foreach \p in {7,...,13}{ \fill (\p ,0) circle (2pt);  \draw (\p,0) circle (4pt);}
\foreach \p in {16,...,22}{ \fill (\p ,0) circle (2pt);  \draw (\p,0) circle (4pt);}
\foreach \p in {1,...,5}{ \fill (\p ,1) circle (2pt);}
\foreach \p in {8,...,14}{ \fill (\p ,1) circle (2pt);}
\foreach \p in {17,...,23}{ \fill (\p ,1) circle (2pt);}
\foreach \p in {2,3,4,9,10,11,12,18,19,20,21}{ \fill (\p ,2) circle (2pt);}
\foreach \p in {0,...,4}{ \fill (\p ,3.5) circle (2pt);  \draw (\p,3.5) circle (4pt);}
\foreach \p in {8,...,12}{ \fill (\p ,3.5) circle (2pt);  \draw (\p,3.5) circle (4pt);}
\foreach \p in {17,...,21}{ \fill (\p ,3.5) circle (2pt);  \draw (\p,3.5) circle (4pt);}
\foreach \p in {1,...,5}{ \fill (\p ,4.5) circle (2pt);}
\foreach \p in {9,...,13}{ \fill (\p ,4.5) circle (2pt);}
\foreach \p in {18,...,22}{ \fill (\p ,4.5) circle (2pt);}
\foreach \p in {2,...,5}{ \fill (\p ,5.5) circle (2pt);}
\foreach \p in {10,...,13}{ \fill (\p ,5.5) circle (2pt);}
\foreach \p in {19,...,21}{ \fill (\p ,5.5) circle (2pt);}
\draw (0,3.5)--(2,5.5)--(3,4.5)--(2,3.5)--(1,3.5)--(3,5.5)--(5,5.5)--(4,4.5)--(5,4.5)--(4,3.5)--(3,3.5);
\draw (8,3.5)--(10,5.5)--(11,4.5)--(10,3.5)--(9,3.5)--(10,4.5)--(11,3.5)--(12,3.5)--(13,4.5)--(12,4.5)--(13,5.5)--(11,5.5);
\draw (17,3.5)--(19,5.5)--(20,4.5)--(21,5.5)--(20,5.5)--(21,4.5)--(22,4.5)--(21,3.5)--(20,3.5)--(19,4.5)--(18,3.5)--(19,3.5);
\draw (0,0)--(1,1)--(2,0)--(1,0)--(2,1)--(3,0)--(4,0)--(5,1); \draw (3,1)--(2,2); \draw [thick, dotted]  (3,1)--(3,2)--(4,2)--(5,1)--(3,1);
\draw (7,0)--(9,2)--(10,1)--(9,0)--(8,0)--(10,2)--(12,2)--(11,1)--(12,0)--(13,0)--(14,1)--(12,1)--(11,0)--(10,0);
\draw (16,0)--(18,2)--(19,1)--(18,0)--(17,0)--(18,1)--(19,0)--(20,0)--(21,1)--(23,1)--(22,0)--(21,0)--(19,2)--(21,2);
\draw  [-stealth] plot [smooth] coordinates {(4,3.5) (4.3,3.5) (4.6,3.8)};
\draw  [-stealth] plot [smooth] coordinates {(12,3.5) (12.3,3.5) (12.6,3.8)};
\draw  [-stealth] plot [smooth] coordinates {(21,3.5) (21.3,3.5) (21.6,3.8)};
\draw  [-stealth] plot [smooth] coordinates {(4,0) (4.3,0) (4.6,0.3)};
\draw  [-stealth] plot [smooth] coordinates {(13,0) (13.3,0) (13.6,0.3)};
\draw  [-stealth] plot [smooth] coordinates {(22,0) (22.3,0) (22.6,0.3)};
\foreach \p in {12,21} {\draw [very thick] (\p,1)--(\p+1,1) ; \draw [-stealth]  (\p+0.5,1) -- (\p+0.8, 1.3) ; }
\end{tikzpicture}
\end{center}
\end{figure}
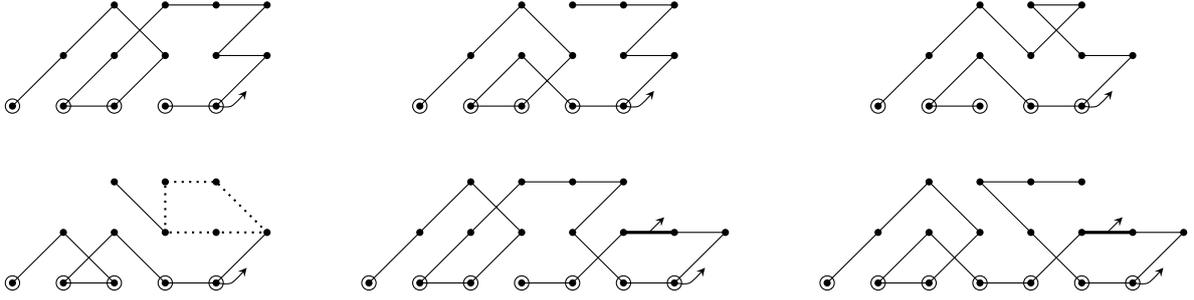

Finally, the constructions needed for $y \in \{5,7\}$ are given in Figure~\ref{fig:smb_small2}.  These work in a similar way to the main constructions.
\end{proof}

\begin{lem}\label{lem:small_c} {\rm ($k=2$; small $c$)}
Let $L = \{1^y, (y-2)^b, y^c\}$ with~$c \leq y-2$ and~$b+c \geq y$.  Then~$L$ has a standard linear realization when~$a \geq y$. 
\end{lem}

\begin{proof}
The framework is similar to the previous lemma.  If either~$b$ or~$c$ is~0 we may use the $\omega$-constructions of Section~\ref{sec:background}, so assume that $b, c > 0$.  Using Lemma~\ref{lem:concat1}, it is sufficient to find standard linear realizations for any value of~$a$ with~$a \leq y-2$; in some cases, the constructions have~$a$ strictly less than~$y-2$ giving a slightly stronger result when the relevant conditions apply.
Similarly, sometimes the constructions naturally cover cases with~$c > y-2$.

Write~$b+c = q(y-2) + r$ with $0 \leq r < y-2$.  The conditions imply that~$q \geq 1$.   The proof splits into four cases depending on the parities of~$y$ and~$r$.

{\em Case 1: $y-2$ and~$r$ even.} 
For all~$c$ in the range $0 \leq c \leq y-2$ we can construct the required standard linear realization of type~$\mathcal{C}_{y-2}$ with~$b+c=y-2$ as a direct application of Lemma~\ref{lem:switch}: the path $\Theta_{y-1,y,b,c}$ realizes~$\{1^{y-1}, (y-2)^b, y^c\}$, is Hamiltonian and is of type $\mathcal{C}_{y-2}$ in~$K_{2y-2}$.  

The first diagram of Figure~\ref{fig:smc1} illustrates this for~$y-2 = 8$; elongate or shorten the trapezoid for the general picture.

\begin{figure}
\caption{Standard linear realizations of type $\mathcal{C}_{y-2}$ for {\em Cases~1} and~{\em 3} of Lemma~\ref{lem:small_c}. }\label{fig:smc1}
\begin{center}
\begin{tikzpicture}[scale=0.8, every node/.style={transform shape}]
\foreach \p in {0,...,9}{ \fill (\p ,1) circle (2pt);}
\foreach \p in {11,...,19}{ \fill (\p ,1) circle (2pt);}
\foreach \p in {1,...,8}{ \fill (\p ,2) circle (2pt);}
\foreach \p in {12,...,18}{ \fill (\p ,2) circle (2pt);}
\draw [thick, dotted] (0,1)--(9,1)--(8,2)--(1,2)--(0,1) ; 
\draw [thick, dotted] (11,1)--(18,1)--(17,2)--(12,2)--(11,1) ; 
\draw (18,1) -- (19,1) -- (18,2);
\foreach \p in {1,3,...,7} {\draw [very thick] (\p,2)--(\p+1,2) ; \draw [-stealth]  (\p+0.5,2) -- (\p+0.2, 2.3) ;}
 \draw [-stealth]  (18,2) -- (17.7, 2.3) ;
 \foreach \p in {12,14,16} {\draw [very thick] (\p,2)--(\p+1,2) ; \draw [-stealth]  (\p+0.5,2) -- (\p+0.2, 2.3) ;}
\end{tikzpicture}
\end{center}
\end{figure}
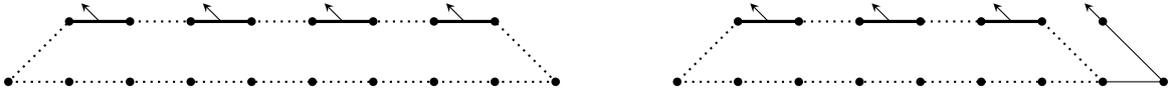

{\em Case 2: $y-2$ even and~$r$ odd.} 
In the {\em Case~1} construction, the final element is~$y-1$.  As in {\em Case~2} of Lemma~\ref{lem:small_b}, we may translate by~1, add an edge from~0 to the final vertex (which is now~$y$) and reverse.  This creates a standard linear realization of type~$\mathcal{C}_{y-2}$ for~$c$ in the range~$1 \leq c \leq y-1$.

{\em Case 3: $y-2$ odd and~$r$ even.} 
For all~$c$ in the range $0 \leq c \leq y-3$ we construct a type~$\mathcal{C}_{y-2}$ standard linear realization with~$b+c=y-2$.  Take the path~$\Theta_{y-2,y,b-1,c}$ in~$K_{2y-2}$ from Lemma~\ref{lem:switch}.  This realizes $\{1^{y-2}, (y-2)^{b-1}, y^c \}$, has final vertex~$y-2$ and omits vertices~$\{ y-1, 2y-3\}$.
Consider the Hamiltonian path $(\Theta_{y-2,y,b-1,c}, y-1, 2y-3)$.  This is a standard linear realization of
$$\{1^{y-2}, (y-2)^{b-1}, y^c \} \cup \{1,y-2\} = 
\{1^{y-1}, (y-2)^{b}, y^c \}$$
and, by the properties in Lemma~\ref{lem:switch} and having final vertex~$2y-3 = v-1$, is of type $\mathcal{C}_{y-2}$.

The second diagram of Figure~\ref{fig:smc1} illustrates this for~$y-2 = 7$; elongate or shorten the trapezoid for the general picture.

{\em Case 4a: $y-2$ even and $r=1$.} 
We construct a $(y-2)$-growable standard linear realization with $b+c = (y-2)+1 = y-1$.  There are two constructions, overlapping in their coverage.

First, assume that~$1 \leq c \leq y-4$.  We use a translation of the path $\Theta'_{y-3,y-1,b-2,c}$ from Lemma~\ref{lem:switch2}, which realizes~$\{1^{y-3}, (y-2)^{b-2}, y^c \}$.  The  Hamiltonian path
$$(y-3, \Theta'_{y-3,y-1,b-2,c}+1, 0)$$
realizes
$$\{ y-2 \} \cup \{1^{y-3}, (y-2)^{b-2}, y^c \} \cup \{y-2\} = \{1^{y-3}, (y-2)^{b}, y^c \}.$$
Reverse it to get the required standard linear realization.  It is $(y-2)$-growable at $2y-4$.  This is illustrated in the first diagram of Figure~\ref{fig:smc_4a}.

Next, take $2 \leq c \leq y-1$.  This time we use a translation of $\Theta_{y-2,y,b,c-2}$ from Lemma~\ref{lem:switch}, which realizes $\{1^{y-2}, (y-2)^b, y^{c-2}\}$.  The Hamiltonian path
$$(y+1, 1, \Theta_{y-2,y,b,c-2} + 2, 0)$$
realizes
$$\{ 1, y \} \cup\{1^{y-2}, (y-2)^b, y^{c-2}\} \cup \{y\} = \{1^{y-1}, (y-2)^{b}, y^c \}.$$
Reverse it to get the required standard linear realization.  It is $(y-2)$-growable at $2y-4$.  This is illustrated in the second diagram of Figure~\ref{fig:smc_4a}.

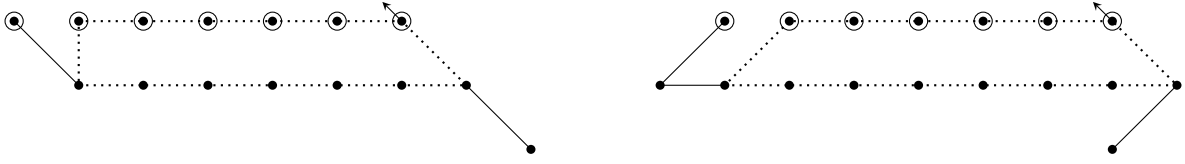
\begin{figure}
\caption{$(y-2)$-Growable standard linear realizations for {\em Case~4a} of Lemma~\ref{lem:small_c}. }\label{fig:smc_4a}
\begin{center}
\begin{tikzpicture}[scale=0.85, every node/.style={transform shape}]
\foreach \p in {1,...,7}{ \fill (\p ,1) circle (2pt);}
\foreach \p in {0,...,6}{ \fill (\p ,2) circle (2pt);}
\foreach \p in {10,...,18}{ \fill (\p ,1) circle (2pt);}
\foreach \p in {11,...,17}{ \fill (\p ,2) circle (2pt);}
\fill (8,0) circle (2pt); \fill (17,0) circle (2pt);
\draw [thick, dotted] (1,1)--(7,1)--(6,2)--(1,2)--(1,1); 
\draw [thick, dotted] (11,1)--(18,1)--(17,2)--(12,2)--(11,1) ;  
\draw (0,2)--(1,1) ; \draw (7,1)--(8,0) ;
\draw (11,2)--(10,1)--(11,1) ; \draw (18,1)--(17,0);
\foreach \p in {0,...,6}{\draw (\p,2) circle (4pt);}
\foreach \p in {11,...,17}{\draw (\p,2) circle (4pt);}
\draw [-stealth] (6,2)--(5.7,2.3);
\draw [-stealth] (17,2)--(16.7,2.3);
\end{tikzpicture}
\end{center}
\end{figure}

{\em Case 4b: $y-2$ and $r$ odd with~$r>1$.} 
The constructions for this subcase is similar to those of Case~4b of Lemma~\ref{lem:small_b}: we construct $(y-2)$-growable standard linear realizations for~$b+c=(y-2)+3 = y-1$, which has $r=3$ and is the smallest target.  These constructions both have edges of the form $(y+i,y+i+1)$ to replace with a three edge sequence adding two more $(y-2)$-edges to cover other values of~$r$.

First take~$1 \leq c \leq y-5$.  We use a translation of the path $\Theta'_{y-5,y-1,b-6,c}$, which realizes~$\{1^{y-5}, (y-2)^{b-6}, y^c\}$.  Consider the Hamiltonian path
$$(y,2,1,y-1,2y-3,2y-2,2y-1,y+1,\Theta'_{y-5,y-1,b-6,c}+3,0),$$
which realizes
$$\{1^3, (y-2)^5\} \cup \{1^{y-5}, (y-2)^{b-6}, y^c\} \cup \{y-2\} = \{1^{y-2}, (y-2)^{b}, y^c\}.$$
The reverse is standard and is $(y-2)$-growable at~$2y-4$.

The construction has 1-edges of the form $(y+i,y+i+1)$  for each even~$i$ with $2 \leq i \leq y-3$.  We may replace some of these with the sequence $(y+i, 2y+i-2, 2y+i-1, y+i+1)$, starting with the lowest value of~$i$ and working up.  Each replacement made adds two $(y-2)$-edges to the linear realization, allowing us to cover the remaining values of~$r$ for these values of~$c$.

Second, take~$4 \leq c \leq y-1$.  We use a translation of the path $\Theta_{y-4,y,b-2,c-4}$, which realizes~$\{1^{y-4}, (y-2)^{b-2}, y^{c-4}\}$.  Consider the Hamiltonian path
$$(y+2,2,1,y+1,2y-1,2y,2y+1,y+3,3 ,\Theta_{y-4,y,b-2,c-4}+4,0),$$
which realizes
$$\{1^3, (y-2)^2, y^3\} \cup \{1^{y-4}, (y-2)^{b-2}, y^{c-4}\} \cup \{y\} = \{1^{y-1}, (y-2)^{b}, y^c\}.$$
The reverse is standard and is $(y-2)$-growable at~$2y-4$.

The construction has 1-edges of the form $(y+i,y+i+1)$  for each even~$i$ with $4 \leq i \leq y-1$.  As in the last construction, we may replace some of these with the sequence $(y+i, 2y+i-2, 2y+i-1, y+i+1)$, starting with the lowest value of~$i$ and working up.  Each replacement made adds two $(y-2)$-edges to the linear realization, allowing us to cover the remaining values of~$r$ for these values of~$c$.

\begin{figure}
\caption{$(y-2)$-Growable standard linear realizations for {\em Case~4b} of Lemma~\ref{lem:small_c}. }\label{fig:smc_4b}
\begin{center}
\begin{tikzpicture}[scale=0.9, every node/.style={transform shape}]
\foreach \p in {1,...,11}{ \fill (\p ,1) circle (2pt);}
\foreach \p in {2,...,10}{ \fill (\p ,2) circle (2pt);}
\foreach \p in {1,...,3}{ \fill (\p ,3) circle (2pt);}
\foreach \p in {3,...,11}{ \fill (\p ,5) circle (2pt);}
\foreach \p in {2,...,10}{ \fill (\p ,6) circle (2pt);}
\foreach \p in {1,...,3}{ \fill (\p ,7) circle (2pt);}
 \fill (10 ,0) circle (2pt);  \fill (12 ,4) circle (2pt);
\foreach \p in {5} {\draw [thick, dotted] (5,\p)--(11,\p)--(10,\p+1)--(5,\p+1)--(5,\p); }
\draw [thick, dotted] (4,1)--(11,1)--(10,2)--(5,2)--(4,1);
\draw (3,2)--(2,1)--(1,1)--(2,2)--(1,3)--(3,3)--(4,2)--(3,1)--(4,1);
\draw (3,6)--(4,5)--(3,5)--(1,7)--(3,7)--(5,5);
\draw (10,0)--(11,1); \draw (11,5)--(12,4);
 \foreach \p in {5,7,9} {\draw [very thick] (\p,6)--(\p+1,6) ; \draw [-stealth]  (\p+0.5,6) -- (\p+0.2, 6.3) ;}
 \foreach \p in {5,7,9} {\draw [very thick] (\p,2)--(\p+1,2) ; \draw [-stealth]  (\p+0.5,2) -- (\p+0.2, 2.3) ;}
\foreach \p in {2,...,10}{\draw (\p,2) circle (4pt);\draw (\p,6) circle (4pt);}
\draw [-stealth] (10,2)--(9.7,2.3);
\draw [-stealth] (10,6)--(9.7,6.3);
\end{tikzpicture}
\end{center}
\end{figure}
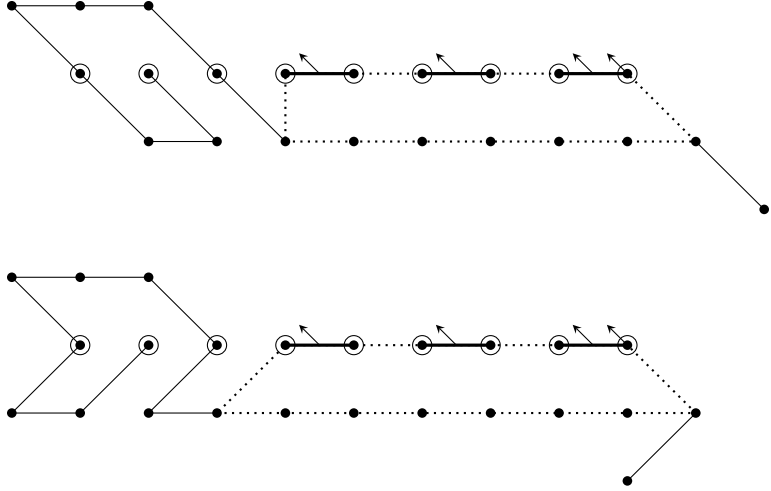

Figure~\ref{fig:smc_4b} illustrates the two constructions, with the same conventions as Figure~\ref{fig:smb_4b}. 
\end{proof}

Bringing these two lemmas together, along with the perfect linear realizations from Section~\ref{sec:para}, we can prove the main result of the section.

\begin{thm}\label{th:y-2} {\rm ($k=2$; $a \geq y$)}
Let $L = \{ 1^a, (y-2)^b, y^c \}$.  $L$ has a linear realization whenever~$a \geq y$.  If $y$ is odd, then~$L$ has a standard linear realization.
\end{thm}

\begin{proof}
First, let~$y$ be odd.
If~$b+c < y$ then we have the required standard linear realization by Lemma~\ref{lem:small_bc}.  Otherwise, if $b \leq y-1$ or $c \leq y-2$ we have the required standard linear realization by Lemma~\ref{lem:small_b} or~\ref{lem:small_c} respectively.  Hence we have to address multisets with~$b > y-1$ and $c > y-2$. 

In this case write
$$L = \{ (y-2)^{ j(y-1) } , y^{ j(y-1) } \} \cup  \{ 1^a, (y-2)^{b'} , y^{c' } \}$$
for the largest integer~$j$ for which it is possible.  Necessarily, either~$b'< y-1$ or~$c' < y-1$ (or both). 

Construct a standard linear realization for~$L$ by concatenating~$j$ perfect linear realizations for $\{ (y-2)^{ (y-1) } , y^{ (y-1) } \}$ from Corollary~\ref{cor:2perf} (valid, as~$y-2$ and~$y$ are coprime) and then concatenating with a standard linear realization for  $\{ 1^a, (y-2)^{b'} , y^{c' } \}$ from one of~Lemma~\ref{lem:small_bc},~\ref{lem:small_b} or~\ref{lem:small_c}.

Now let~$y$ be even.  There are two adjustments needed to the above argument, one beneficial and one detrimental.

On the negative side, we do not have a perfect linear realization from Corollary~\ref{cor:2perf} available, as~$y-2$ and~$y$ are not coprime.  However, we may use the standard linear realization for~$\{1, (y-2)^{j(y-2)}, y^{jy}\}$ from Theorem~\ref{th:tracks_div} in its place.

On the positive side, the constructions of Lemmas~\ref{lem:small_bc},~\ref{lem:small_b} and~\ref{lem:small_c} in fact show that we can construct a standard linear realization for~$\{ 1^a, (y-2)^{b'} , y^{c' } \}$, where~$b'<y-2$ and~$c' < y$, provided that~$a \geq y-1$.

The same argument as for odd~$y$ now gives the stated result.
\end{proof}

As previously noted, for even~$y$ we have the result from~\cite{OPPS} that $\{ 1^a, (y-2)^b, y^c \}$ has a standard linear realization for~$a \geq y-1$.  Hence Theorem~\ref{th:y-2} is sufficient to prove Theorem~\ref{th:biggie}.  The alternative constructions that combine to give the version in Theorem~\ref{th:y-2} add to the possibilities for future generalizations to larger~$k$.

\section{Concluding remarks on future directions}\label{sec:future}

The main aim of this section is re-emphasize the potential of our techniques in context and highlight opportunities for future inquiries that could inspire further work.

Our current investigations show promise of a systematic approach to generalize our constructions for $k=1$ and $k=2$ to arbitrary odd and even $k$, respectively (maybe even beyond the restrained scope of Conjecture~\ref{conj:y-k}), with a similar bound on $a$. This would allow us to cover the majority of the cases of the Coprime BHR Conjecture for supports of size 3.

There is also a line of inquiry available to reduce the bound on $a$ for a given $k$. We saw in the introduction that Theorem~\ref{th:biggie} can complete previously open instances of the Coprime BHR Conjecture for some supports and values of~$v$.   As a first step towards a more thorough understanding of when this happens---and to get a sense of what supports future work should investigate---we consider the equivalent supports to~$\{1, y-k, y\}$.

\begin{lem}\label{lem:equiv} {\rm (Equivalent supports)}
Let~$U = \{1, y-k, y\}$ be a strongly admissible support of size~$v-1$.  Then the two equivalent supports are of the form~$\{ 1, \widehat{s}, \widehat{ks - 1} \}$ and $\{ 1, \widehat{t}, \widehat{kt + 1} \}$ for~$s = y^{-1}$ and~$t = (y-k)^{-1}$.
\end{lem}

\begin{proof}
Consider the supports~$U'$ and~$U''$ obtained by multiplying by~$y^{-1}$ and by~$(y-k)^{-1}$ respectively, and then reducing: 
$$U' = \{ \widehat{y^{-1}}, \widehat{y^{-1}(y-k)}, 1 \} 
  = \{ 1, \widehat{y^{-1}}, \widehat{ky^{-1}-1} \}$$ 
$$U'' = \{ \widehat{(y-k)^{-1}}, 1,  \widehat{(y-k)^{-1}y} \} 
  = \{ 1, \widehat{(y-k)^{-1}}, \widehat{k(y-k)^{-1}+1} \}$$ 
as required.
\end{proof}

Hence we are motivated to consider supports of the form~$\{1, \widehat{z}, \widehat{kz \pm 1}\}$ for small~$k$.  It is often the case that $\widehat{z} = z$ and $\widehat{kz \pm 1} = kz \pm 1$.  Therefore, the supports~$\{1,z,kz\pm 1\}$ are natural first ones to consider.  We believe that the techniques developed in the present paper may be adapted to this situation and that aiming for a result of the form ``multisets of the form $\{1^a, z^b, (kz\pm 1)^c\}$ are realizable for all~$a \geq a^*$," for some fixed~$a^* \approx kz$, is a promising line of future enquiry.

Should such a line of enquiry be successful, work along the lines of that done in Section~\ref{sec:1y-1y} here will probably then be required.  Specifically, some result similar to Lemma~\ref{lem:sawtooth} to provide standard linear realizations for $\{1^a, (y-k)^b, y^c\}$ when~$b$ and~$c$ are small for~$a$ as small as possible.  While the ``sawtooth" pattern used there does not immediately generalize without complications, we point to ways in which other ideas in the paper might be applied to the problem.  For brevity, we limit the comments to~$k=2$ and do not formalize the results.

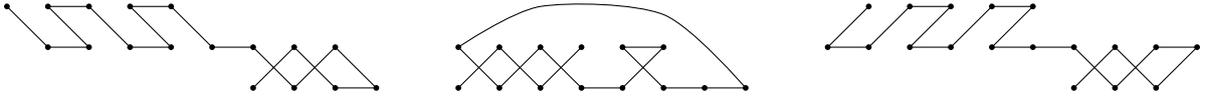
\begin{figure}
\caption{Construction diagrams for $\{1^a, (y - 2)^b, y^c\}$ with $b + c$ small.}\label{fig:k2smbc}
\begin{center}
\begin{tikzpicture}[scale=0.54, every node/.style={transform shape}]
\foreach \p in {6,...,9}{ \fill (\p ,0) circle (2pt);}
\foreach \p in {11,...,18}{ \fill (\p ,0) circle (2pt);}
\foreach \p in {26,...,28}{ \fill (\p ,0) circle (2pt);}
\foreach \p in {1,...,8}{ \fill (\p ,1) circle (2pt);}
\foreach \p in {11,...,16}{ \fill (\p ,1) circle (2pt);}
\foreach \p in {20,...,29}{ \fill (\p ,1) circle (2pt);}
\foreach \p in {0,...,4}{ \fill (\p ,2) circle (2pt);}
\foreach \p in {21,...,25}{ \fill (\p ,2) circle (2pt);}
\draw (0,2)--(1,1)--(2,1)--(1,2)--(2,2)--(3,1)--(4,1)--(3,2)--(4,2)--(5,1)--(6,1)--(7,0)--(8,1)--(9,0)--(8,0)--(7,1)--(6,0);
\draw (11,0)--(12,1)--(13,0)--(14,1); \draw (11,1)--(12,0)--(13,1)--(14,0)--(15,0)--(16,1)--(15,1)--(16,0)--(18,0); 
\draw   plot [smooth] coordinates {(11,1) (13,2) (16,1.8) (18,0)};
\draw (21,2)--(20,1)--(21,1)--(22,2)--(23,2)--(22,1)--(23,1)--(24,2)--(25,2)--(24,1)--(26,1)--(27,0)--(28,1)--(29,1)--(28,0)--(27,1)--(26,0);
\end{tikzpicture}
\end{center}
\end{figure}

Figure~\ref{fig:k2smbc} gives illustrative constructions for the cases when~$b+c$ is even.  The first diagram corresponds to~$b>c$, the second to~$b=c$ (the rightmost ``$\times$" shape is not needed when~$b$ and~$c$ are odd), and the third to~$b < c$.  When~$b+c$ is even, the diagrams obtained by removing the final two edges from each of the first and third diagram give the necessary approaches (no analog of the second is required).

Using these construction methods, it is possible to obtain standard linear realizations for~$\{1^a, (y-2)^b, y^c\}$ with~$b+c<y$ with a lower bound for~$a$ of $y - \min(b,c)$.

\end{document}